\documentclass[10pt]{elsarticle}
\usepackage{lineno,hyperref}
\usepackage{booktabs}
\usepackage{multirow}
\usepackage{amsthm}
\usepackage{graphicx}
\usepackage{subfigure}
\usepackage{float}
\usepackage{bm}
\usepackage{lipsum}
\usepackage{listings} 
\usepackage{color}
\usepackage{amsfonts,enumerate,amsmath,amssymb}
\def\qed{\strut\hfill $\Box$}
\newtheorem{thm}{Theorem}[section]

\newtheorem{lem}[thm]{Lemma}

\newcommand{\thmref}[1]{Theorem~{\rm \ref{#1}}}
\newcommand{\lemref}[1]{Lemma~{\rm \ref{#1}}}

\def\para#1{\vskip .4\baselineskip\noindent{\bf #1}}
\numberwithin{equation}{section}
\begin{document}
	\begin{frontmatter}

		\title{Large deviation principle for a two-time-scale McKean-Vlasov model with  jumps}
		
		\author[mymainaddress]{Xiaoyu Yang}
		\ead{yangxiaoyu@yahoo.com}
		
		\author[mymainaddress]{Yong Xu\corref{mycorrespondingauthor}}
		\cortext[mycorrespondingauthor]{Corresponding author}
		\ead{hsux3@nwpu.edu.cn}

		\address[mymainaddress]{School of Mathematics and Statistics, Northwestern Polytechnical University, Xi'an, 710072, China}

		\begin{abstract}
This work focus on the  large deviation principle for a two-time scale McKean-Vlasov system with jumps.  Based on the variational framework of the McKean-Vlasov system with jumps, it is turned into weak convergence for the controlled system.  Unlike  general two-time scale system,  the controlled McKean-Vlasov system is related to the law of the original system, which causes difficulties in qualitative analysis.  In solving this problem, employing asymptotics of the original system and a Khasminskii-type averaging principle together is efficient.  Finally, it is shown that the limit is related to the Dirac measure of the solution to the ordinary differential equation.
			
			\vskip 0.08in			
		\end{abstract}		
	\end{frontmatter}

	\section{Introduction}\label{sec-1}
	The McKean-Vlasov system can be traced back to the original work of stochastic toy model related to the Vlasov kinetic system of plasma by Kac \cite{1956Foundations}. Shortly afterwards, McKean researched the propagation of chaos in  interacting particle system, which is related to Boltzmann's model for the statistical mechanics of rarefied gases \cite{1967Propagation}. 
	Take the number of particles go to infinity, then the above particle systems converge to the mean-field system, which is the well-known  McKean-Vlasov system. {The McKean-Vlasov system  does not only depend on the solution itself but also depend on its time marginal law.} Up to now, the McKean-Vlasov system has attracted a lot of attention since it has widely employed  in several fields, including biology, physics, chemistry, and so on. As regards  properties of  the solution to such systems, see for instance \cite{2020Mishura,2021Rockner,2018Wang}.
	
	Furthermore, an enormous number of problems in physics and mechanics can be reduced to two-time scale systems, which consist of two or more subsystems with different time scales. Consequently, this work focuses on the two-time scale Mckean-Vlasov system with jumps as follows,
	\begin{eqnarray}\label{1}
	\left
	\{
	\begin{array}{ll}
	dX^{\varepsilon ,\delta}_t = b_1(X^{\varepsilon ,\delta}_t,\mathcal{L}_{X^{\varepsilon ,\delta}_t}, Y^{\varepsilon ,\delta}_t)dt + \sqrt \varepsilon  \sigma_{1}( X^{\varepsilon ,\delta}_t,\mathcal{L}_{X^{\varepsilon ,\delta}_t})dW_{t}+ \varepsilon \int_{\mathbf{X}} g(t, X^{\varepsilon ,\delta}_t ,\mathcal{L}_{X^{\varepsilon ,\delta}_t},z)\tilde{N}^{\frac{1}{\varepsilon }}(dzdt),\\
	 dY^{\varepsilon ,\delta}_t =
	 \frac{1}{\delta} b_2( X^{\varepsilon ,\delta}_t,\mathcal{L}_{X^{\varepsilon ,\delta}_t}, Y^{\varepsilon ,\delta}_t)dt + \frac{1}{{\sqrt \delta}} \sigma_2( X^{\varepsilon ,\delta}_t, \mathcal{L}_{X^{\varepsilon ,\delta}_t}, Y^{\varepsilon ,\delta}_t)dW_{t},
	\end{array}
	\right.
	\end{eqnarray}
where $t\in [0,T]$, $(X^{\varepsilon ,\delta}_0, Y^{\varepsilon ,\delta}_0)=(X_0, Y_0)\in \mathbb{R}^{n}\times \mathbb{R}^{n}$ and $t\in [0,T]$. $W$ is a $\mathbb{R}^d-$valued Brownian motion (Bm). Independent of Bm $W$, $\tilde{N}^{\frac{1}{\varepsilon }}(dzdt)=N^{\frac{1}{\varepsilon }}(dzdt)-\frac{1}{\varepsilon}\nu(dz) dt$ is the compensated Poisson random measure with
	   associated Poisson measure $N^{\frac{1}{\varepsilon }}(dzdt)$, intensity measure $\frac{1}{\varepsilon}\nu(dz) dt$, L\'evy measure  $\nu$  satisfying $\int_\mathbf{X}(1\wedge z^2)\nu(dz)< \infty$ in a locally compact Polish space $\mathbf{X}$ \cite{Applebaum}. 
	 {$\mathcal{L}_{X^{\varepsilon ,\delta}_t}$ stands for the distribution of slow variable $\{X^{\varepsilon ,\delta}_t\}$ for $t\in [0,T]$}.	$\{X^{\varepsilon,\delta}\}$ is called the slow component and $\{Y^{\varepsilon,\delta}\}$ is the fast component.  $\varepsilon$ and $\delta$ are small  parameters satisfying $0<\delta=o(\varepsilon)<1$, which are used to describe the separation of different time scales.  For $\mu \in \mathcal{P}$ where  $\mathcal{P}$ is the set of all probability measure on $(\mathbb{R}^d, \mathcal{B}(\mathbb{R}^d))$, Set
	 $$\mathcal{P}_{2}:=\big\{\mu \in \mathcal{P}: \mu\left(|\cdot|^{2}\right):=\int_{\mathbb{R}^{n}}|x|^{2} \mu(d x)<\infty\big\},$$
	 Then the set $\mathcal{P}_{2}$ is a Polish space under the $L^2$-Wasserstein distance, 
	 $$\mathbb{W}_{2}\left(\mu_{1}, \mu_{2}\right):=\inf _{\pi \in \mathcal{C}_{\mu_{1}, \mu_{2}}}\left[\int_{\mathbb{R}^{n} \times \mathbb{R}^{n}}|x-y|^{2} \pi(d x, d y)\right]^{1 / 2},\quad \text{for}\quad \mu_{1}, \mu_{2} \in \mathcal{P}_{2}$$
	 where $\mathcal{C}_{\mu_{1}, \mu_{2}}$ is the set of all couplings of measures $\mu_1$ and $\mu_2$, i.e. $\pi\in \mathcal{C}_{\mu_{1}, \mu_{2}}$ is a probability measure on $\mathbb{R}^{n} \times \mathbb{R}^{n}$ satisfying that $\pi(\cdot \times \mathbb{R}^{n})=\mu_{1}$ and $\pi(\mathbb{R}^{n} \times \cdot  )=\mu_{2}$.
	 Then, 
	 $b_{i}:  \mathbb{R}^{n}\times\mathcal{P}_2\times \mathbb{R}^{n} \rightarrow  \mathbb{R}^{n}$, $\sigma_{1}: \mathbb{R}^{n}\times\mathcal{P}_2 \rightarrow \mathbb{R}^{ n\times n}$, $\sigma_{2}: \mathbb{R}^{n}\times\mathcal{P}_2\times \mathbb{R}^{n} \rightarrow \mathbb{R}^{ n\times n}$, $g: [0,T] \times\mathbb{R}^{n}\times\mathcal{P}_2\times \mathbf{X} \rightarrow  \mathbb{R}^{n}$ are nonlinear functions.

The large deviation is an important topic in the field of probability \cite{Dembo2009Large}.  As a complement and development of the Law of Large Numbers (LLN) and the Central Limit Theorem (CLT), large deviation principles could characterise the exponential decay rate of rare event probabilities \cite{2010Asymptotic}. Moreover, large deviations have wide applications in statistics, complex systems engineering and so on \cite{2009The}. 
Large deviation principles for stochastic dynamical systems under small noise were proposed by Freidlin and Wentzell \cite{1984Random}. Subsequently, the large deviation principle has been studied intensively, see \cite{1992Large,2004Large,2006Large,2009Large} and the references are given there.	
Up to now, there have been several kinds of methods to study the large deviation principle, such as the weak convergence method \cite{1997A,2014Large,2008Large,DS,2011Variational,2013Large,2020Large}, the PDE theory \cite{BCS}, the nonlinear semigroup theory and the viscosity solution approach proposed in \cite{KP,FK, FFK}. For the McKean-Vlasov system driven by standard Bm, by the exponential equivalence arguments, the distribution in the original McKean-Vlasov system can be replaced by the Dirac measure of the solution of ordinary differential equations (ODE) \cite{Dos2019Freidlin,Liu2021Long}. For the McKean-Vlasov model with jumps, however, it is difficult to find similar results. Fortunately, by constructing the variational framework for the above system, then the weak convergence method could be constructed for the large deviation principle of the McKean-Vlasov model \cite{W2022Large}. Then, along this weak convergence approach, large deviation principles for McKean-Vlasov quasilinear stochastic evolution systems were established \cite{Hong2021Large}.

	 However, there are just a few works focusing on the two-time-scale McKean-Vlasov system, and all these works just aimed at the Bm \cite{Hong2021Central,Suo2021Central}. Hence, we mainly study  large deviation principles for a two-time scale McKean-Vlasov system with jumps (\ref{1}).
 In our work, based on the weak convergence approach with respect to the McKean-Vlasov system with jumps in \cite{W2022Large}, the problem could turn into basic qualitative properties (in other words, weak convergence) for the controlled system. In detail, it is related to the distribution of the original system, but not of the controlled system. 
Therefore, unlike general two-time scale systems, it is necessary to treat the probability distributions of the original slow component when it comes to analysing the controlled system. Before proving the target result, we can show that the original slow system strongly converges to the averaged ODE. Next, the combination of the strong convergence of the original system and the Khasminskii-type averaging principle is used to efficiently analyse the controlled McKean-Vlasov model.    With the particular regime that $\delta=o(\varepsilon)$ we could see that in the weak limit there is no control in the fast component. The weak convergence of the controlled slow component is obtained by the property of the original system, the Burkholder-Davis-Gundy inequality, It\^o's formula and the exponential ergodicity of the fast component without a controlled term. It is observed that the limit is related to the Dirac measure of the solution of the ODE.

	The paper is organized  as follows. In  Section 2, we set up notations and  some precise conditions for the two-time scale system (\ref{1}), and state  our result. Section \ref{sec-3}  reviews some preliminary results. The proof of our main result in Section 4. 
	Throughout this paper, $c$, $C$, $c_1$, $C_1$, $\cdots$ denote certain positive constants which may vary from line to line. Denote $C([0, T]; \mathbb{R}^{n})=\mathbf{C}$ be the space of continuous functions,  and $\mathbf{D}=D([0, T], \mathbb{R}^{n})$ be the space of $\mathbb{R}^{n}$-valued, c\`{a}dl\`ag functions endowed with  the  Skorohod topology. 
	
	\section{Notations, Assumptions and  Main Results}\label{sec-2}
		\subsection{Preliminaries and Notations}\label{sec-2-1}
	Set  $\mathcal{B}(\mathbf{X})$ be the Borel $\sigma$-field on locally compact Polish space $\mathbf{X}$. Set $\mathcal{M}_{F}(\mathbf{X})$ be the space of all Borel measure $\nu$ on $\mathbf{X}$  with $\nu(K)< \infty$ for  compact subset $K\subset\mathbf{X}$. $\mathcal{M}_{F}(\mathbf{X})$ is the Polish space under the topology that $\langle f, v\rangle=\int_{\mathbf{X}} f(u) \nu(d u), \nu \in \mathcal{M}_{F}(\mathbf{X})$ for every $f\in C_c(\mathbf{X})$ (continuous function space with compact support).
	 
	Let $\lambda_{T}$ be Lebesgue measure on $[0, T]$.  
	Let  $\mathbf{M}=\mathcal{M}_{F}([0, T]\times \mathbf{X})$, then denote  $\mathbb{P}$ the  probability measure on $(\mathbf{M}, \mathcal{B}(\mathbf{M}))$    under  the Poisson random measure $N(m)=m: \mathbf{M}\rightarrow\mathbf{M}$ with intensity measure $\lambda_{T}\otimes\nu$. 
	We denote the product space $\mathbf{V}=\mathbf{C}\times\mathbf{M}$. Let $W=(w_{i})_{i=1}^{d}: w_{i}(w, m)=w_{i}$ be coordinate maps on $\mathbf{V}$. Define the Poisson random measure $N(w, m)=m: \mathbf{V}\rightarrow\mathbf{M}$. Now, $W$ is independent of Poisson random measure $N$.	 Set $ \mathcal{H}_t=\sigma\{N((0, s]\times A), \beta_{i}: s\leq t, A\in\mathcal{B}(\mathbf{X})\}$.

 Let $\bar{\mathbf{M}}=\mathcal{M}_{F}([0, T]\times \mathbf{X}\times \mathbb{R}_{+})$, then $\bar{\mathbb{P}}$ is the  probability measure on $(\bar{\mathbf{M}}, \mathcal{B}(\bar{\mathbf{M}}))$ with Poisson random measure $\bar{N}(m)=m: \bar{\mathbf{M}}\rightarrow \bar{\mathbf{M}}$  with intensity  $\bar{\nu}_{T}=\lambda_{T}\otimes\nu\otimes\lambda_{\infty}$ where $\lambda_{\infty}$ is Lebesgue measure on $\mathbb{R}_{+}$. 
	Let $\bar{\mathbf{V}}=\mathbf{C}\times\bar{\mathbf{M}}$, then we  define the  Poisson random measure $\bar{N}$ and Brownian motion  $W=(w_{i})_{i=1}^{d}$ on $\bar{\mathbf{V}}$ analogously. Further, set $(\bar{\mathbb{P}}, \bar{\mathcal{H}_t})$ on $(\bar{\mathbf{V}}, \mathcal{B}(\bar{\mathbf{V}}))$.
	Here and subsequently, denote by $\bar{\mathcal{F}}_{t}$ the   $\bar{\mathbb{P}}$-completion of the filtration $\bar{\mathcal{H}_t}$, and   $\bar{\mathcal{P}}$  the predictable $\sigma$-field on $[0, T]\times\bar{\mathbf{V}}$ with the filtration $\{\bar{\mathcal{F}}_{t}: 0\leq t\leq T\}$ on $(\bar{\mathbf{V}}, \mathcal{B}(\bar{\mathbf{V}}))$. 
	
	Set $U=(U^{i})_{i=1,...,d}\in L^2([0, T]; \mathbb{R}^{n})$ with  norm
	\[
	\int_{0}^{T}\|U(s)\|^2ds=\int_{0}^{T}\big(\sum_{i=1}^{d}|U^{i}(s)|^{2}\big)ds <\infty, \quad \textrm{a.s. $\bar{\mathbb{P}}$.}
	\]
	For each $U\in L^2([0, T]; \mathbb{R}^{n})$, set $L^{(1)}(U)=\frac{1}{2}\int_{0}^{T}\|U(s)\|^2ds$.

 Set $\ell(r)=r\log r-r+1: [0, \infty)\rightarrow [0, \infty)$.	Let $\bar{\mathcal{A}}$ be the class of all $\bar{\mathcal{P}}\otimes\mathcal{B}(\mathbf{X})-\mathcal{B}[0, \infty)$ {measurable} maps $V: [0, T]\times\bar{\mathbf{V}}\times \mathbf{X}\rightarrow [0, \infty)$. Since $(\bar{\mathbf{V}}, \mathcal{B}(\bar{\mathbf{V}}))$ is underlying probability space, we will replace $V(t, w, m, z)$, $(w, m)\in\bar{\mathbf{V}}$ by $V(t, z)$ for simplicity.
		 For each $V\in\bar{\mathcal{A}}$, define $L^{(2)}(V)$ by
	\begin{eqnarray*}
		L^{(2)}(V)(\omega)=\int_{[0,T]\times\mathbf{X}}\ell(V(t, z))\nu_{T}(dtdz).
	\end{eqnarray*}

Set $\mathcal{U}= L^2([0, T]; \mathbb{R}^{n})\times \bar{\mathcal{A}}$. 	For each $(U, V)\in\mathcal{U}$, denote that
	\[
	L(U,V)=L^{(1)}(U) + L^{(2)}(V).
	\] 
	For $m\in\mathbb{N}$, let 
	\[
	S_{1}^{m}=\{U\in L^2([0, T]; \mathbb{R}^{n}): L^{(1)}(U)\leq m \},
	\]
	and
	\[
	S_{2}^{m}=\{V\in \bar{\mathcal{A}}: L^{(2)}(V)\leq m \}.
	\]
 Let $\mathbf{S}=\bigcup_{m\in \mathbb{N}}\big(S_{1}^{m}\times S_{2}^{m} \big)$ and $\mathcal{U}^{m}$ be the space of controls, that is
	\[
	\mathcal{U}^{m}=\{(U, V)\in\mathcal{U}: (U, V)\in S_{1}^{m}\times S_{2}^{m}, \bar{\mathbb{P}}-a.e.\}.
	\]
	
\subsection{Assumptions and Main Results}\label{sec-2-2}
	
We give assumptions needed in next section.
	\begin{itemize}
		\item[\textbf{A1}.]  There exists a constant $C_1> 0$ such that for  any $ (x_1,\mu_1, y_1) $,  $ (x_2, \mu_2,y_2)\in \mathbb{R}^{n}\times\mathcal{P}_{2} \times\mathbb{R}^{n}$, 
		\begin{eqnarray*}
		\begin{aligned}
			&|b_1( x_2,\mu_2,y _2) - b_1( x_1,\mu_1,y _1)|^2 +|b_2( x_2,\mu_2, y_2) - b_2( x_1,\mu_1,y _1)|^2+|\sigma_1( x_2,\mu_2)-\sigma_1(x_1,\mu_1)|^2\\
		&+|\sigma_2( x_2, \mu_2,y_2)-\sigma_2( x_1,\mu_1, y_1)|^2+\int_{\mathbf{X}}| g(t, x_2,\mu_2, z)-g(t, x_1,\mu_1, z)|^2\nu(dz)\\
		&\leq C_1(|x_2-x_1|^2+|y_2-y_1|^2+\mathbb{W}^2_2(\mu_1,\mu_2)).
		\end{aligned}
		\end{eqnarray*}	
	\end{itemize}	
Due to Assumption 	(\textbf{A1}), it could deduce that 
there exists a constant $ C_2>  0$ such that for  all $(x, \mu,y)  \in \mathbb{R}^{n}\times\mathcal{P}_{2} \times\mathbb{R}^{n}$, 
\begin{eqnarray*}
	|b_1( x,\mu, y)|^{2}+|b_2( x,\mu, y)|^{2}+|\sigma_1( x,\mu)|^2+\int_{\mathbf{X}}| g(t, x,\mu, z)|^2\nu(dz)\leq C_2 (1+|x|^2+|y|^2+\mu(|\cdot|^2)),
\end{eqnarray*}
holds.

Under Assumption (\textbf{A1}), for  initial value $(X^{\varepsilon ,\delta}_0, Y^{\varepsilon ,\delta}_0)=(X_0, Y_0)\in \mathbb{R}^{n}\times \mathbb{R}^{n}$, there exists a unique strong solution $(X^{\varepsilon ,\delta},Y^{\varepsilon ,\delta})$ in $\mathbf{D}\times \mathbf{C}$  to  the two-time scale McKean-Vlasov system (\ref{1}), which is from \cite[Chapter 6]{Applebaum}. Then there exists a measurable map 
	\[
	\mathcal{G}^{\varepsilon ,\delta}(\sqrt \varepsilon W, \varepsilon N^{\frac{1}{\varepsilon }}): \mathbf{C}\times\mathbf{M}\rightarrow \mathbf{D}
	\]
	such that
	$X^{\varepsilon ,\delta}:=\mathcal{G}^{\varepsilon,\delta}(\sqrt \varepsilon W, \varepsilon N^{\frac{1}{\varepsilon }})$. 
	
Moreover, here follow other assumptions.	
\begin{itemize}
	\item[\textbf{A2}.] There exists a constant $ C_3>  0$ such that for  all $(x, \mu,y)  \in \mathbb{R}^{n}\times\mathcal{P}_{2} \times\mathbb{R}^{n}$,
	\begin{eqnarray*}
		\sup_{y\in\mathbb{R}^{n}}|\sigma_2(x,\mu, y)|^2\leq C_3 (1+|x|^2+\mu(|\cdot|^2)),.
	\end{eqnarray*}		
holds.
		\item[\textbf{A3}.]  There exists a constant $,C_4,C_5,C_6 > 0$ such that for any  $(x, \mu,y_1),(x, \mu,y_2),(x, \mu,y)   \in \mathbb{R}^{n}\times\mathcal{P}_{2} \times\mathbb{R}^{n}$
		\[
			2\left\langle y _1 - y _2, b_2( x, \mu,y _1) - b_2( x,\mu, y _2) \right\rangle + |\sigma_2(x,\mu, y_2)-\sigma_2( x,\mu, y_1)|^2 \le  - C_4 \vert  {{y _1} - {y _2}}  \vert ^2,
		\]
and
		\[
			\left\langle {y ,{b_2}\left(  x,\mu, y \right)} \right\rangle + |\sigma_2( x, \mu,y)|^2 \le  - C_5{ \vert  y   \vert ^2}+C_6(1+|x|^2+\mu(|\cdot|^2)),
		\]
		hold.
		\item[\textbf{A4}.]  There exists a constant $\varrho\in(0, \infty)$ such that for all $E \in \mathcal{B}([0, T]\times \mathbf{X})$, $\nu_{T}(E)<\infty$
		\[
			\int_{E}e^{\varrho\|g(t, z)\|} \nu_{T}(dzdt) < \infty
		\]
		with $\|g(t, z)\|=\{\sup_{x\in\mathbb{R}^{n},\mu\in \mathcal{P}(\mathbb{R}^d)}\frac{|g(t, x,\mu, z)|^2}{1+|x|^2+\mu(|\cdot|^2)}\}$.
	\end{itemize}		
According to similar arguments to \cite[Lemma 3.6, Propositon 3.7]{M2021Strong},  Assumption (\textbf{A3}) could ensure that the solution to the following fast equation with frozen-{$(X,\mu)$}, 
\[d\tilde{Y}_t = b_2({X}, \mu,\tilde{Y}_t) dt +\sigma_2( {X},\mu, \tilde{Y}_t)dW_{t}\]
has a unique  invariant probability measure $\mu_{{X}}$. Moreover,
	$\mathcal{L}_{\bar X_t}=\delta_{\bar X_t}$ is the {Dirac} measure for the solution to the following ODEs,
	\begin{eqnarray}\label{1-1}
	d{\bar X}_t = \bar{b}_1(\bar{X}_t,\delta_{ \bar X_t})dt,
	\end{eqnarray}
	with $\bar  X_0=X_0$ and 	 {$\bar{b}_1(\cdot)=\int_{\mathbb{R}^{n}}b_1(\cdot, {\tilde Y})\mu_{\cdot}(d{\tilde Y})$}. 
	 For any $X_0\in \mathbb{R}^{n}$, there exists a unique solution ${\bar X}$ to the above deterministic ODE  \eqref{1-1}.
	 
Then we could define the skeleton equation as follows
	\begin{eqnarray}\label{3}
	d\hat{X}_t = \bar{b}_1(\hat{X}_t,\mathcal{L}_{\bar X_t})dt + \sigma_{1}( \hat{X}_t,\mathcal{L}_{\bar X_t})\psi_tdt+ \int_{\mathbf{X}} g(t, \hat{X}_t,\mathcal{L}_{\bar X_t}, z)(\phi_t-1)\nu(dz)dt.
	\end{eqnarray}
	From deterministic equation (\ref{3}) we could  define the solution map
	\[
	\mathcal{G}^{0}: S_{1}^{m}\times S_{2}^{m} \rightarrow C([0,T]; \mathbb{R}^{n})
	\]
	such that $\hat{X}=\mathcal{G}^{0}(\psi, \phi)$.

		Now, the statement of  main theorem is given.	
	\begin{thm}\label{thm}
		  Assume {(\textbf{A1})}--{(\textbf{A4})},  $\delta =o (\varepsilon)$, we let $\varepsilon \to 0$.
		 The slow variable $X^{\varepsilon ,\delta}$ of two-time scale McKean-Vlasov model (\ref{1}) satisfies the large deviation principle on $\mathbf{D}$ with the good rate function $I: \mathbf{D}\rightarrow [0, \infty)$
		\begin{eqnarray}\label{rate}
		I(\xi) = \inf_{(\psi,\phi)\in S_{\xi }} L(\psi,\phi),
		\end{eqnarray}
		where $S_{\xi }:=\{(\psi,\phi)\in \mathbf{S}: \xi=\mathcal{G}^{0}(\psi,\phi)\}$ for $\xi\in\mathbf{D}$.
	\end{thm}	
	The proof of \thmref{thm} will be shown in Section 4.
	
	\section{Preliminary Lemmas}\label{sec-3}

	Before proving  \thmref{thm},  we  give some prior estimates.

		\begin{lem}\label{lem1}
			Under  Assumptions {(\textbf{A1})}--{(\textbf{A3})}, for any $(X_0, Y_0)\in \mathbb{R}^{n}\times \mathbb{R}^{n}$, and $t\in [0,T]$, we have
			\begin{eqnarray}\label{1--1}
           \mathbb{E}\big[ {\mathop {\sup }\limits_{0 \le t \le T} | {X}_{t}^{\varepsilon ,\delta} |^2} \big]
          <\infty ,\qquad
           \mathbb{E}\big[ { | {Y}_{t}^{\varepsilon ,\delta} |^2} \big]<\infty.
			\end{eqnarray}
		\end{lem}
		\para{Proof}. {According to  the It\^o's formula, we can get 
		\begin{eqnarray}\label{lemma2.10}
		\begin{aligned}
		\mathbb{E}[{| {Y}_{t}^{\varepsilon ,\delta} |^2}] &=\mathbb{E}{[| Y_0 |^2] } + \frac{2}{\delta }\mathbb{E}\int_0^t {\langle {Y}_{s}^{\varepsilon ,\delta},b_2( {X}_{s}^{\varepsilon ,\delta},\mathcal{L}_{X^{\varepsilon ,\delta}_s},{Y}_{s}^{\varepsilon ,\delta} )\rangle ds}   +\frac{2}{\sqrt{\delta} }\mathbb{E}\int_0^t {\langle {X}_{s}^{\varepsilon ,\delta},\sigma_{2}( {{X}_{s}^{\varepsilon ,\delta},\mathcal{L}_{X^{\varepsilon ,\delta}_s},{Y}_{s}^{\varepsilon ,\delta}} ) dW_s}\rangle\\
		&\quad+ \frac{1}{\delta }\mathbb{E}\int_0^t |{{\sigma_{2}}( {{X}_{s}^{\varepsilon ,\delta},\mathcal{L}_{X^{\varepsilon ,\delta}_s},{Y}_{s}^{\varepsilon ,\delta}} )|^2ds}.
		\end{aligned}
		\end{eqnarray}}	
	where $|\sigma_2|$ is the Hilbert-Schmidt norm of the matrix $\sigma_2$.
	
It is easy to see that the fourth term is a true martingale. Then, we have $\mathbb{E}[\int_0^t {\langle {X}_{s}^{\varepsilon ,\delta},\sigma_{2}( {{X}_{s}^{\varepsilon ,\delta},\mathcal{L}_{X^{\varepsilon ,\delta}_s},{Y}_{s}^{\varepsilon ,\delta}} ) dW_s\rangle}]=0$.
	Then, we have
	\begin{eqnarray*}\label{3-21}
	\begin{aligned}
	\frac{d\mathbb{E}[{|  {Y}_{t}^{\varepsilon ,\delta} |^2}]}{dt} &= \frac{2}{\delta }\mathbb{E} {\langle  {Y}_{t}^{\varepsilon ,\delta},b_2(  {X}_{t}^{\varepsilon ,\delta},\mathcal{L}_{X^{\varepsilon ,\delta}_t}, {Y}_{t}^{\varepsilon ,\delta} )\rangle } + \frac{1}{\delta }\mathbb{E} |{\sigma_{2} }( { {X}_{t}^{\varepsilon ,\delta},\mathcal{L}_{X^{\varepsilon ,\delta}_t}, {Y}_{t}^{\varepsilon ,\delta}} )|^2.
	\end{aligned}
	\end{eqnarray*} 
	With Assumption (\textbf{A4}), we have
	\begin{eqnarray*}\label{lemma2.12}
	&\frac{2}{\delta }{\langle  {Y}_{t}^{\varepsilon ,\delta},b_2(  {X}_{t}^{\varepsilon ,\delta},\mathcal{L}_{X^{\varepsilon ,\delta}_t}, {Y}_{t}^{\varepsilon ,\delta} )\rangle } +\frac{1}{\delta } |{\sigma_{2} }( { {X}_{t}^{\varepsilon ,\delta},\mathcal{L}_{X^{\varepsilon ,\delta}_t}, {Y}_{t}^{\varepsilon ,\delta}} )|^2\cr
	&\le	- \frac{{ 2C_5 }}{\delta } {{| {Y}_{t}^{\varepsilon ,\delta} |^2}}  + \frac{{ C_6}}{\delta }(1+{| {X}_{t}^{\varepsilon ,\delta} |^2}+\mathcal{L}_{X^{\varepsilon ,\delta}_t}(|\cdot|^2))	.	
	\end{eqnarray*}
	Thus,  we have
		\begin{eqnarray*}\label{lemma2.16}
	\begin{aligned}
	\frac{d\mathbb{E}[{|  {Y}_{t}^{\varepsilon ,\delta} |^2}]}{dt}  &\le 		- \frac{{ 2C_5 }}{\delta } \mathbb{E}[{{| {Y}_{t}^{\varepsilon ,\delta} |^2}}]  + \frac{{ C_6}}{\delta }(1+\mathbb{E}[{| {X}_{t}^{\varepsilon ,\delta} |^2}]+\mathcal{L}_{X^{\varepsilon ,\delta}_t}(|\cdot|^2))	.	
	\end{aligned}
	\end{eqnarray*}
	Moreover,  by comparison theorem, we have for all $t$ that
	\begin{eqnarray}\label{lemma2.17}
	\begin{aligned}
	\mathbb{E}[{|  {Y}_{t}^{\varepsilon ,\delta} |^2}] &\le  |y_0|^2  e^{-\frac{2C_5}{\delta} t}+  \frac{{ C_6}}{\delta }\int_{0}^{t}e^{-\frac{2C_5(t-s)}{\delta}}  (1+\mathbb{E}[{| {X}_{s}^{\varepsilon ,\delta} |^2}]+\mathcal{L}_{X^{\varepsilon ,\delta}_s}(|\cdot|^2))ds.
	\end{aligned}
	\end{eqnarray}
		After taking the expectation on the both sides of  (\ref{lemma2.17}), and by the  Gronwall's inequality, it leads to that 
		\begin{eqnarray}\label{lemma2.171}
		\mathbb{E}[{| {Y}_{t}^{\varepsilon ,\delta} |^2}] &\le& {c_1}\mathbb{E}\mathop {\sup }\limits_{s \in \left[ {0,t} \right]} {| {X}_{s}^{\varepsilon ,\delta} |^2}+c_2.
		\end{eqnarray}
	By  the It\^o's formula, we get
		\begin{eqnarray}\label{lemma2.211}
		{| {X}_{t}^{\varepsilon ,\delta} |^2}={| {x_0} |^2}  +{\mathcal{G}_1} + {\mathcal{G}_2} + {\mathcal{G}_3} + {\mathcal{G}_{4}}+ {\mathcal{G}_{5}},
		\end{eqnarray}
		where 
		\begin{eqnarray*}
			\begin{aligned}
				&{\mathcal{G}_1}= 2\int_0^t {\langle {X}_{s}^{\varepsilon ,\delta},b_1( {{X}_{s}^{\varepsilon ,\delta},\mathcal{L}_{X^{\varepsilon ,\delta}_s},{Y}_{s}^{\varepsilon ,\delta}} )\rangle ds},\\
				&
				{\mathcal{G}_2}= 2\sqrt \varepsilon  \int_0^t {\langle {X}_{s}^{\varepsilon ,\delta},\sigma_{1}( {X}_{s}^{\varepsilon ,\delta},\mathcal{L}_{X^{\varepsilon ,\delta}_s} )dW_s\rangle } ,\\
				&
				{\mathcal{G}_3}= \varepsilon \int_0^t {|{\sigma_{1}}( {X}_{s}^{\varepsilon ,\delta},\mathcal{L}_{X^{\varepsilon ,\delta}_s})|^2ds},\\
				&{\mathcal{G}_{4}}= \int_0^t {\int_{\mathbf{X}} {\big[ {{( {X}_{s}^{\varepsilon ,\delta}{\rm{ + }}\varepsilon g(s, {X}_{s}^{\varepsilon ,\delta} ,\mathcal{L}_{X^{\varepsilon ,\delta}_s},z) )^2} - {| {X}_{s}^{\varepsilon ,\delta} |^2}} \big]{{ \tilde N}^{{{1} \mathord{\left/{\vphantom {{\phi^{\varepsilon ,\delta} \left( s \right)} \varepsilon }} \right.\kern-\nulldelimiterspace} \varepsilon }}}\left( {dzds} \right)} },\\
				&{\mathcal{G}_{5}}=\varepsilon\int_0^t {\int_{\mathbf{X}} {|g(s, {X}_{s}^{\varepsilon ,\delta},\mathcal{L}_{X^{\varepsilon ,\delta}_s},z )|^2 \nu(dz)ds} }.
			\end{aligned}
		\end{eqnarray*}
		By conditions (\textbf{A1}), we get that
		\begin{eqnarray}\label{lemma2.3}
		\begin{aligned}
		\left| {{\mathcal{G}_1}} \right| &\le \left( {{C_2} + 1} \right)\int_0^t {| {X}_{s}^{\varepsilon ,\delta} |^2ds}  + {C_2}T + {C_2}\int_0^t \mathcal{L}_{X^{\varepsilon ,\delta}_s}(|\cdot|^2)ds+ {C_2}\int_0^t {| {Y}_{s}^{\varepsilon ,\delta} |^2ds},\\	
		\left| {{\mathcal{G}_3}} \right| &\le 2\varepsilon C_2 \int_0^t { {| {X}_{s}^{\varepsilon ,\delta} |^2}ds}  + \varepsilon {C_2}T,\\
		\left| {{\mathcal{G}_{5}}} \right| &\le  \varepsilon {C_2} {\int_0^t { ( 1 + | {X}_{s}^{\varepsilon ,\delta} |^2 +\mathcal{L}_{X^{\varepsilon ,\delta}_s}(|\cdot|^2))ds } }.
		\end{aligned}
		\end{eqnarray}
		Estimates (\ref{lemma2.211}) and (\ref{lemma2.3}) yield that
		\begin{eqnarray*}
			\begin{aligned}
				{| {X}_{t}^{\varepsilon ,\delta} |^2} &\le {| x_0 |^2}+ {C_2}T +\left( 2{C_2} + 1 +\varepsilon \right)\int_0^t { | {X}_{s}^{\varepsilon ,\delta} |^2ds}  + \left( {1 + {C_2}} \right)\int_0^t \big[ 1 + | {X}_{s}^{\varepsilon ,\delta} |^2 +\mathcal{L}_{X^{\varepsilon ,\delta}_s}(|\cdot|^2)\big]ds\\
				&\quad+ {C_2}\int_0^t {| {Y}_{s}^{\varepsilon ,\delta} |^2ds} + | {{\mathcal{G}_{2}}} |+ | {{\mathcal{G}_{4}}} |.
			\end{aligned}
		\end{eqnarray*}
		Then with aid of the Gronwall's lemma, we can conclude that
		\begin{eqnarray*} \label{lemma2.212}
		\mathop {\sup }\limits_{t \in \left[ {0,T} \right]}{| {X}_{t}^{\varepsilon ,\delta} |^2} &\le& c_3\big({| x_0|^2}+\mathop {\sup }\limits_{t \in \left[ {0,T} \right]} | {{\mathcal{G}_2}} |+\mathop {\sup }\limits_{t \in \left[ {0,T} \right]} | {{\mathcal{G}_{4}}} |+ {C_2}T + {C_2}\int_0^T {| {Y}_{s}^{\varepsilon ,\delta} |^2ds}\big).
		\end{eqnarray*}		
		Note that  the term ${{\rm \mathcal{G}}_4}$ can be rearranged as follows
		\begin{eqnarray*}\label{lemma2.51}
		{{\rm \mathcal{G}}_4} := {{\rm \mathcal{G}}_{41}} + {{\rm \mathcal{G}}_{42}},
		\end{eqnarray*}
		where
		\begin{eqnarray*}\label{lemma2.52}
		\begin{aligned}
		{{\rm \mathcal{G}}_{41}} 
		&= {\varepsilon ^2}\int_0^t {\int_{\mathbf{X}} {  |{g(s,{X}_{s}^{\varepsilon ,\delta},\mathcal{L}_{X^{\varepsilon ,\delta}_s} ,z )}|^2  {\tilde N^{{1 \mathord{\left/{\vphantom {{\phi^{\varepsilon ,\delta} } \varepsilon }} \right.\kern-\nulldelimiterspace} \varepsilon }}}\left( {dzds} \right)} },\\
		{{\rm \mathcal{G}}_{42}}
		&=2\varepsilon \int_0^t {\int_{\mathbf{X}} \langle{ {{\rm \mathcal{G}}_s^{\varepsilon ,\delta }},g(s, {X}_{s}^{\varepsilon ,\delta} ,\mathcal{L}_{X^{\varepsilon ,\delta}_s}, z)}\rangle {\tilde N}^{1 \mathord{\left/{\vphantom {{\phi^{\varepsilon ,\delta} } \varepsilon }} \right.\kern-\nulldelimiterspace} \varepsilon }( {dzds} ) },
		\end{aligned}
		\end{eqnarray*}		
		For the term ${{\rm \mathcal{G}}_{41}}$, by Assumptions (\textbf {A2}) and (\textbf {A5}), we have
		\begin{eqnarray*}
		\begin{aligned}
		\mathbb{E}[\mathop {\sup }\limits_{t \in [ {0, T} ]}  {{{\rm \mathcal{G}}_{41}}} ] &\le \mathbb{E}\Big[\mathop {\sup }\limits_{t \in [ {0, T} ]}  \varepsilon \int_0^t \int_{\mathbf{X}} |g(s, {X}_{s}^{\varepsilon ,\delta},\mathcal{L}_{X^{\varepsilon ,\delta}_s}, z )|^2  \nu(dz)ds \Big]\\	
		&\le \varepsilon \mathbb{E}\Big[  \int_0^{ T}  ( 1 + | {X}_{s}^{\varepsilon ,\delta} |^2 +\mathcal{L}_{X^{\varepsilon ,\delta}_s}(|\cdot|^2))ds   \Big] \\
		&\le \varepsilon c_{3}  .
		\end{aligned}	
		\end{eqnarray*}
	By using the Burkholder-Davis-Gundy inequality, we get
		\begin{eqnarray*}
		\begin{aligned}	
		\mathbb{E}[\mathop {\sup }\limits_{t \in[ {0, T} ]} {{{\rm \mathcal{G}}_{42}}} ] &\le 4\mathbb{E}\left[ {{\rm \mathcal{G}}_{42}^2} \right]_{T}^{{1 \mathord{\left/{\vphantom {1 2}} \right.\kern-\nulldelimiterspace} 2}}\\		
		&\le 8\mathbb{E}{\Big[ {\varepsilon \mathop {\sup }\limits_{t \in [ {0, T} ]} | {{\rm \mathcal{G}}_t^{\varepsilon ,\delta }} |^2\int_0^{ T} {\int_{\mathbf{X}} {{g^2}(s, {X}_{s}^{\varepsilon ,\delta},\mathcal{L}_{X^{\varepsilon ,\delta}_s}, z ) \nu( dz )ds} } } \Big]^{{1 \mathord{\left/{\vphantom {1 2}} \right.\kern-\nulldelimiterspace} 2}}}\\
		&\le \frac{1}{{8c_{3}}}\mathbb{E}[\mathop {\sup }\limits_{t \in [ {0, T} ]} {| {{\rm \mathcal{G}}_t^{\varepsilon ,\delta }}|^2}] + c_{3}\varepsilon \mathbb{E}[\sup_{t\in[0,T]}|X_t^{\varepsilon ,\delta}|^2] + c_{3}\varepsilon  T.
		\end{aligned}
		\end{eqnarray*}		
		Due to  Burkholder-Davis-Gundy inequality, it deduces that
		\begin{eqnarray}\label{lemma2.7}
		\begin{aligned}
		\mathbb{E}[\mathop {\sup }\limits_{t \in \left[ {0,T} \right]} \left| {{\mathcal{G}_{4}}} \right|] 
		&\le \varepsilon {c_3}\mathbb{E}[{\mathop {\sup }\limits_{t \in \left[ {0,T} \right]} | {X}_{t}^{\varepsilon ,\delta} |^2}]  +\varepsilon {c_3}.
		\end{aligned}
		\end{eqnarray}
		Using  the Gronwall's lemma, from the above it follows  that
		\begin{eqnarray*}\label{lemma2.102}
			\mathbb{E}[\mathop {\sup }\limits_{t \in \left[ {0,T} \right]} {| {X}_{t}^{\varepsilon ,\delta} |^2}] \le c_4.
		\end{eqnarray*}		
		This proof is completed.\qed

		\begin{lem}\label{lem2}
			Under  Assumptions {(\textbf{A1})}--{(\textbf{A4})} and $t\in [0,T]$, for $t(\Delta):={\left[ \frac{t}{ \Delta}\right] }\Delta$, we have
			\begin{eqnarray}\label{1-2}
            \mathbb{E}\big[ {|{X}_{t}^{\varepsilon ,\delta} - { {X}_{t(\Delta)}^{\varepsilon ,\delta}}|^2}\big]<C\Delta (1+|x_0|+|y_0|).
			\end{eqnarray}
		\end{lem}
			\para{Proof}. 	
				Indeed,  it has
				\begin{eqnarray*}\label{lemma2.21}
				{{X}_{t}^{\varepsilon ,\delta} - { {X}_{t(\Delta)}^{\varepsilon ,\delta}}}
				&=& \int_{t\left( \Delta  \right)}^t b_1(X^{\varepsilon ,\delta}_s,\mathcal{L}_{X^{\varepsilon ,\delta}_s}, Y^{\varepsilon ,\delta}_s)ds + \int_{t\left( \Delta  \right)}^t \sqrt \varepsilon  \sigma_{1}( X^{\varepsilon ,\delta}_s,\mathcal{L}_{X^{\varepsilon ,\delta}_s})dW_{s}\cr
				&&+\varepsilon\int_{t\left( \Delta  \right)}^t \int_{\mathbf{X}} g(s, X^{\varepsilon ,\delta}_s ,\mathcal{L}_{X^{\varepsilon ,\delta}_s},z)\tilde{N}^{\frac{1}{\varepsilon }}(dzds)\cr
				&=:&{{\mathcal A}_1} + {{\mathcal A}_2} + {{\mathcal A}_3} .
				\end{eqnarray*}
				Using Assumption {(\textbf{A1})} and  H\"older inequality, it follows that
				\begin{eqnarray*}\label{lemma2.23}
				\mathbb{E}[|{{\mathcal A}_1} |^2]
				&\le&C_2 \Delta\int_{t\left( \Delta  \right)}^t {( 1 + | {X}_{s}^{\varepsilon ,\delta} |^2 +\mathcal{L}_{X^{\varepsilon ,\delta}_s}(|\cdot|^2)+|  {Y}_{s}^{\varepsilon ,\delta} |^2)ds}.
			\end{eqnarray*}
			Then, by the It\^o isometry, we get
		\begin{eqnarray*}\label{lemma2.23a}			
				\mathbb{E}[|{{\mathcal A}_2}|^2]   &\le& 2\varepsilon {C_2}\mathbb{E}\int_{t\left( \Delta  \right)}^t {\big( 1 + |  {X}_{s}^{\varepsilon ,\delta} |^2+\mathcal{L}_{X^{\varepsilon ,\delta}_s}(|\cdot|^2) \big)ds} .
			\end{eqnarray*}
		With Assumption (\textbf{A1}) and  Burkholder-Davis-Gundy inequality, we could see that
			\begin{eqnarray*}\label{lemma2.23b}
				\mathbb{E}[|{{\mathcal A}_3}|^2] &\le& \varepsilon C_2\int_{t\left( \Delta  \right)}^t { {\big( 1 + |  {X}_{s}^{\varepsilon ,\delta} |^2+\mathcal{L}_{X^{\varepsilon ,\delta}_s}(|\cdot|^2) \big)}ds}.
				\end{eqnarray*}			
				Thus, from what has already been proved in \lemref{lem1}, it deduces  that (\ref{1-2}) holds.

			The proof is completed.\qed

		\begin{lem}\label{lem3}
		Under  Assumptions {(\textbf{A1})}--{(\textbf{A3})}, and let $\varepsilon\to 0$.  {The slow variable $X^{\varepsilon ,\delta}$ of original two-time scale McKean-Vlasov system (\ref{1}) strongly converges to $\bar x$, which is the solution to the ODE (\ref{1-1}) as $\varepsilon \to  0$, i.e.}
		\begin{eqnarray}\label{1-3}
		\lim_{\varepsilon,\delta\to 0}\mathbb{E}[\sup_{t\in[0,T]}|X^{\varepsilon,\delta}_t -\bar X_t|^2] =0.
		\end{eqnarray}
	   \end{lem}
		\para{Proof}.
Before proving \eqref{1-3}, we construct  the auxiliary processes as follows, for $t(\Delta):={\left[ \frac{t}{ \Delta}\right] }\Delta$,
\begin{eqnarray*}
	\left
	\{
	\begin{array}{ll}
		d\check{X}_{t}^{\varepsilon ,\delta} &=  b_1( {{ {X}_{t(\Delta)}^{\varepsilon ,\delta}},\mathcal{L}_{X^{\varepsilon ,\delta}_{t(\Delta)}},\check{Y}_t^{\varepsilon ,\delta }} )dt, \\
		d\check{y}_t^{\varepsilon ,\delta } &=   \frac{1}{\delta }b_2\big( {{ {X}_{t(\Delta)}^{\varepsilon ,\delta}},\mathcal{L}_{X^{\varepsilon ,\delta}_{t(\Delta)}},\check{Y}_t^{\varepsilon ,\delta }} \big)dt + \frac{1}{{\sqrt \delta  }}\sigma_{2} \big( {{ {X}_{t(\Delta)}^{\varepsilon ,\delta}},\mathcal{L}_{X^{\varepsilon ,\delta}_{t(\Delta)}},\check{Y}_t^{\varepsilon ,\delta }} \big)dW_t .
	\end{array}
	\right.
\end{eqnarray*}
Take  similar arguments in \lemref{lem1}, it has
\begin{eqnarray*}
	\mathbb{E}\big[ {\mathop {\sup }\limits_{0 \le t \le T} | \check{X}_{t}^{\varepsilon ,\delta} |^2} \big]< \infty,\qquad
	\mathbb{E}\big[ {| {\check{Y}_t^{\varepsilon ,\delta }} |^2} \big]< \infty.
\end{eqnarray*}
{Next, our task is now to show that $$\mathbb{E}[\sup_{t \in[0, T]}|{ {X}_{t}^{\varepsilon ,\delta}}-{ \check{X}_{t}^{\varepsilon ,\delta}}|^2]\le C\Delta.$$
	Here, $C>0$ is a constant independent of $\varepsilon,\delta,\Delta$. 
Let ${{\rm \mathcal{M}}_t^{\varepsilon ,\delta }}:={ {X}_{t}^{\varepsilon ,\delta}}-{ \check{X}_{t}^{\varepsilon ,\delta}}$.
By the It\^o's formula, it deduces that
\begin{eqnarray*}\label{lemma2.41}
	{| {{\rm \mathcal{M}}_t^{\varepsilon ,\delta }}|}^2= {{\rm \mathcal{M}}_1} + {{\rm\mathcal{M}}_2} + {{\rm \mathcal{M}}_3} + {{\rm \mathcal{M}}_4}  + {{\rm \mathcal{M}}_5},
\end{eqnarray*}} 
where
\begin{eqnarray*}\label{lemma2.80}
	\begin{aligned}
		{{\rm \mathcal{M}}_1} &= 2\int_0^t {\langle {{\rm \mathcal{M}}_s^{\varepsilon ,\delta }} ,[ {b_1( {{X}_{s}^{\varepsilon ,\delta},\mathcal{L}_{X^{\varepsilon ,\delta}_s},{Y}_{s}^{\varepsilon ,\delta}} ) - b_1( {{x}_{s( \Delta  )}^{\varepsilon ,\delta},\mathcal{L}_{X^{\varepsilon ,\delta}_{s(\Delta)}},\check{Y}_s^{\varepsilon ,\delta }} )} ]\rangle ds}\\
		{{\rm \mathcal{M}}_2}&=  \varepsilon \int_0^t {|{\sigma_{1}}( {X}_{s}^{\varepsilon ,\delta},\mathcal{L}_{X^{\varepsilon ,\delta}_s} )|^2ds},\\
		{{\rm \mathcal{M}}_3}
		&= 2\sqrt \varepsilon  \int_0^t {\langle  {{\rm \mathcal{M}}_s^{\varepsilon ,\delta }}, \sigma_{1}( {X}_{s}^{\varepsilon ,\delta},\mathcal{L}_{X^{\varepsilon ,\delta}_s} ) dW_s}\rangle, \\
		{{\rm \mathcal{M}}_4}
		&= \int_0^t \int_{\mathbf{X}} {[ {( {{{\rm \mathcal{M}}_s^{\varepsilon ,\delta }}{\rm{ + }}\varepsilon g(s, {X}_{s}^{\varepsilon ,\delta} ,\mathcal{L}_{X^{\varepsilon ,\delta}_s},z)} )^2} - ( {{\rm \mathcal{M}}_s^{\varepsilon ,\delta }} )^2 ]{\tilde N^{\frac{1}{\varepsilon}}( dzds )} },\\
		{{\rm \mathcal{M}}_5}
		&= \varepsilon\int_0^t {\int_{\mathbf{X}}  |g(s, {X}_{s}^{\varepsilon ,\delta} ,\mathcal{L}_{X^{\varepsilon ,\delta}_s},z)|^2 \nu(dz)ds }, 
	\end{aligned}
\end{eqnarray*}
Assumption {(\textbf{A1})} and elementary inequality yield that
\begin{eqnarray}\label{lemma2.42}
\begin{aligned}
{{\rm \mathcal{M}}_1} 
&\le \int_{0}^t ( {{\rm \mathcal{M}}_s^{\varepsilon ,\delta }} )^2ds  + {{\rm \mathcal{M}}_{11}},\\
{{\rm \mathcal{M}}_2}+{{\rm \mathcal{M}}_5}  
&\le C_2\varepsilon\int_0^t ( 1+|X_s^{\varepsilon ,\delta}|^2+ \mathcal{L}_{X^{\varepsilon ,\delta}_s}(|\cdot|^2))ds,
\end{aligned}
\end{eqnarray}
where
\[
{{\rm \mathcal{M}}_{11}}
:=\int_0^t [ {b_1( {{X}_{s}^{\varepsilon ,\delta},\mathcal{L}_{X^{\varepsilon ,\delta}_s},{Y}_{s}^{\varepsilon ,\delta}} ) - b_1( {{X}_{s( \Delta  )}^{\varepsilon ,\delta},\mathcal{L}_{X^{\varepsilon ,\delta}_{s(\Delta)}},\check{Y}_s^{\varepsilon ,\delta }} )} ]^2ds.
\]
Note that  the term ${{\rm \mathcal{M}}_4}$ can be rearranged as follows
\begin{eqnarray}
	{{\rm \mathcal{M}}_4} := {{\rm \mathcal{M}}_{41}} + {{\rm \mathcal{M}}_{42}},
\end{eqnarray}
where
\begin{eqnarray}\label{lemma2.44}
	\begin{aligned}
		{{\rm \mathcal{M}}_{41}} 
		&= {\varepsilon ^2}\int_0^t {\int_{\mathbf{X}} {  {g^2( s,{X}_{s}^{\varepsilon ,\delta},\mathcal{L}_{X^{\varepsilon ,\delta}_s} ,z )}  {\tilde N^{{1 \mathord{\left/{\vphantom {{\phi^{\varepsilon ,\delta} } \varepsilon }} \right.\kern-\nulldelimiterspace} \varepsilon }}}\left( {dzds} \right)} },\\
		{{\rm \mathcal{M}}_{42}}
		&=2\varepsilon \int_0^t {\int_{\mathbf{X}} \langle{ {{\rm \mathcal{M}}_s^{\varepsilon ,\delta }},g(s, {X}_{s}^{\varepsilon ,\delta} ,\mathcal{L}_{X^{\varepsilon ,\delta}_s}, z)}\rangle {\tilde N}^{1 \mathord{\left/{\vphantom {{\phi^{\varepsilon ,\delta} } \varepsilon }} \right.\kern-\nulldelimiterspace} \varepsilon }( {dzds} ) },
	\end{aligned}
\end{eqnarray}
With aid of  the  Gronwall's inequality, it implied  from (\ref{lemma2.42}) to (\ref{lemma2.44}), 
\begin{eqnarray*}\label{lemma2.46}
	{| {X}_{t}^{\varepsilon ,\delta} - \check{X}_{t}^{\varepsilon ,\delta} |^2} \le {e^{\left(   2t + {C_1} \right)}}\left\{ {{{\rm \mathcal{M}}_{11}} + {{\rm \mathcal{M}}_3}  + {{\rm \mathcal{M}}_{41}} + {{\rm \mathcal{M}}_{42}} } \right\},
\end{eqnarray*}
which shows that		
\begin{eqnarray}\label{lemma2.47}
\mathbb{E}[\mathop {\sup }\limits_{t \in [ {0, T} ]} {|  {X}_{t}^{\varepsilon ,\delta} - \check{X}_{t}^{\varepsilon ,\delta}  |^2}] \le c_{7}\mathbb{E}[\mathop {\sup }\limits_{t \in [ {0, T} ]} \left( {{{\rm \mathcal{M}}_{11}} + {{\rm \mathcal{M}}_3} + {{\rm \mathcal{M}}_{41}} + {{\rm \mathcal{M}}_{42}}} \right)],
\end{eqnarray}
which is from choosing  the constant $c_{7}\ge{e^{T}}$.
{By the  \cite[Lemma 3.4]{M2021Strong}, it has that 
	for $\varepsilon, \delta> 0$ small enough,   }
\begin{eqnarray}\label{lemma2.86}
\mathbb{E}\Big[\int_0^{ T}{| {{ {Y}_{t}^{\varepsilon ,\delta}} - \check{Y}_t^{\varepsilon ,\delta }} |^2}dt\Big] \le c_{8}\Delta_1,
\end{eqnarray}
where $\Delta_1$ is small enough and related to $\Delta$.

Next, by Assumption (\textbf{A1}) and estimate (\ref{lemma2.86}), it follows that 
\begin{eqnarray}\label{lemma2.48}
\begin{aligned}
\mathbb{E}[\mathop {\sup }\limits_{t \in [ {0, T} ]}  {{{\rm \mathcal{M}}_{11}}}] &\le\mathbb{E}\Big[\mathop {\sup }\limits_{t \in [ {0, T} ]} \int_0^t [ {b_1( {{X}_{s}^{\varepsilon ,\delta},\mathcal{L}_{X^{\varepsilon ,\delta}_s},{Y}_{s}^{\varepsilon ,\delta}} ) - b_1( {{X}_{s\left( \Delta  \right)}^{\varepsilon ,\delta},\mathcal{L}_{X^{\varepsilon ,\delta}_{s(\Delta)}},\check{Y}_s^{\varepsilon ,\delta }} )} ]^2ds\Big] \\
&\le {C_1}\mathbb{E} \Big[\int_0^{T} \big[ ( {X}_{s}^{\varepsilon ,\delta} - {X}_{s\left( \Delta  \right)}^{\varepsilon ,\delta} )^2 + ({Y}_{s}^{\varepsilon ,\delta} - \check{Y}_s^{\varepsilon ,\delta } ) ^2 +\mathbb{W}^2_2(\mathcal{L}_{X^{\varepsilon ,\delta}_s},\mathcal{L}_{X^{\varepsilon ,\delta}_{s(\Delta)}})\big]ds \Big]\\
&\le c_{9}\Delta_1.
\end{aligned}
\end{eqnarray}
For the term ${{\rm \mathcal{M}}_{41}}$, by Assumptions (\textbf {A2}), we have
\begin{eqnarray}\label{lemma2.50}
\begin{aligned}
\mathbb{E}[\mathop {\sup }\limits_{t \in [ {0, T} ]}  {{{\rm \mathcal{M}}_{41}}} ] &\le \varepsilon\mathbb{E}\Big[\mathop {\sup }\limits_{t \in [ {0, T} ]}   \int_0^t \int_{\mathbf{X}}  g^2(s, {X}_{s}^{\varepsilon ,\delta},\mathcal{L}_{X^{\varepsilon ,\delta}_s}, z )  \nu\left( dz \right)ds \Big]\\	
&\le \varepsilon \mathbb{E}\Big[  \int_0^{ T}  ( 1 + | {X}_{s}^{\varepsilon ,\delta} |^2 +\mathcal{L}_{X^{\varepsilon ,\delta}_s}(|\cdot|^2))ds   \Big] \\
&\le \varepsilon c_{10}  .
\end{aligned}	
\end{eqnarray}
Using the Burkholder-Davis-Gundy inequality and elementary inequality, we get
\begin{eqnarray}\label{lemma2.49}
\begin{aligned}	
\mathbb{E}[\mathop {\sup }\limits_{t \in [ {0, T} ]} {{{\rm \mathcal{M}}_3}}] &\le 8\mathbb{E}{\Big[ {\varepsilon \int_0^{ T} {{( {{\rm \mathcal{M}}_s^{\varepsilon ,\delta }} )^2}|{\sigma_{1}}( {X}_{s}^{\varepsilon ,\delta},\mathcal{L}_{X^{\varepsilon ,\delta}_s} )|^2ds} } \Big]^{{1 \mathord{\left/{\vphantom {1 2}} \right.\kern-\nulldelimiterspace} 2}}}\\
&\le \frac{1}{{8c_{7}}}\mathbb{E}[\mathop {\sup }\limits_{t \in [ {0, T} ]} ( {{\rm \mathcal{M}}_t^{\varepsilon ,\delta }} )^2] + \varepsilon c_{12}(1+\mathbb{E}[\sup_{t\in[0,T]}|X_t^{\varepsilon ,\delta}|^2]),\\
\mathbb{E}[\mathop {\sup }\limits_{t \in[ {0, T} ]} {{{\rm \mathcal{M}}_{42}}} ] &\le 4\mathbb{E}\left[ {{\rm \mathcal{M}}_{42}^2} \right]_{T}^{{1 \mathord{\left/{\vphantom {1 2}} \right.\kern-\nulldelimiterspace} 2}}\\		
&\le 8\mathbb{E}{\Big[ {\varepsilon \mathop {\sup }\limits_{t \in [ {0, T} ]} | {{\rm \mathcal{M}}_t^{\varepsilon ,\delta }} |^2\int_0^{ T} {\int_{\mathbf{X}} {{g^2}(s, {X}_{s}^{\varepsilon ,\delta},\mathcal{L}_{X^{\varepsilon ,\delta}_s}, z ) \nu( dz )ds} } } \Big]^{{1 \mathord{\left/{\vphantom {1 2}} \right.\kern-\nulldelimiterspace} 2}}}\\
&\le \frac{1}{{8c_{7}}}\mathbb{E}[\mathop {\sup }\limits_{t \in [ {0, T} ]} {| {{\rm \mathcal{M}}_t^{\varepsilon ,\delta }}|^2}] + c_{12}\varepsilon \mathbb{E}[\sup_{t\in[0,T]}|X_t^{\varepsilon ,\delta}|^2]\\	
&\quad + c_{12}\varepsilon  T.
\end{aligned}
\end{eqnarray}		
Then by estimates (\ref{lemma2.47}),  (\ref{lemma2.48}), (\ref{lemma2.50}) and (\ref{lemma2.49}), we can get
\begin{eqnarray*}
	\mathbb{E}[\mathop {\sup }\limits_{t \in [ {0, T} ]} {| {X}_{t}^{\varepsilon ,\delta} - \check{X}_{t}^{\varepsilon ,\delta}|^2}] \le c_{11}\Delta_1  + \frac{1}{4}\mathbb{E}[\mathop {\sup }\limits_{t \in [ {0, T} ]} {| {X}_{t}^{\varepsilon ,\delta} - \check{X}_{t}^{\varepsilon ,\delta} |^2}],
\end{eqnarray*}
which implies that
\begin{eqnarray}\label{lemma2.53}
\mathbb{E}[\mathop {\sup }\limits_{t \in [ {0, T} ]} {| {X}_{t}^{\varepsilon ,\delta} - \check{X}_{t}^{\varepsilon ,\delta}|^2}] \le \frac{4}{3}c_{11}\Delta_1.
\end{eqnarray}				
Next,  the rest of the proof runs as   \cite[Theorem 2.2]{Hong2022Strong},  we can show that 
\begin{eqnarray}\label{lemma2.531}
\mathbb{E}[\mathop {\sup }\limits_{t \in [ {0, T} ]} {|\check {X}_{t}^{\varepsilon ,\delta} - \bar{X}_{t}|^2}] \le c_{11}\frac{\delta}{\Delta}.
\end{eqnarray}
Combine \eqref{lemma2.53} and \eqref{lemma2.531}, we have
\begin{eqnarray}\label{lemma2.532}
\mathbb{E}[\mathop {\sup }\limits_{t \in [ {0, T} ]} {| {X}_{t}^{\varepsilon ,\delta} - \bar{X}_{t}|^2}] \le c\frac{\delta}{\Delta},
\end{eqnarray}
where $c>0$ is independent of $\varepsilon,\delta,\Delta$.
Choose suitable $\Delta>0$ such that as $\delta \to 0$, $\Delta$ and $\frac{\delta}{\Delta}$ converge to $0$. Then,
the estimate (\ref{1-3}) can be shown.

		 The proof is completed.\qed

	\section{Proof of the Main Result \thmref{thm}}\label{sec-4}
	\para{Proof of the \thmref{thm}}. 
	
		{From \cite[Theorem 4.4]{W2022Large}, it provides a convenient, sufficient condition to  prove large deviations. In what follows, we will show the verification of (a) in \cite[Theorem 4.4]{W2022Large} in Step 1. The verification of (b) in \cite[Theorem 4.4]{W2022Large} will be shown in Step 2.}

	\textbf{Step 1}. The proof in this step is in deterministic sense.
	
	Let $(\psi^{(j)}, \phi^{(j)})$ and $(\psi, \phi)$ belong to $\big(S_{1}^{m}\times S_{2}^{m} \big)$ such that $(\psi^{(j)}, \phi^{(j)})\rightarrow(\psi, \phi)$ as $j\rightarrow\infty$.
	Assume $\{\hat{X}^{(j)}\}$ is a family of solutions to the skeleton equation (\ref{3}), that is,
	\begin{eqnarray} \label{3-1}
		d\hat{X}^{(j)}_t =  \bar{b}_1( \hat{X}^{(j)}_t,\mathcal{L}_{\bar X_t})dt + \sigma_{1}(\hat{X}^{(j)}_t,\mathcal{L}_{\bar X_t})\psi_t^{(j)}dt+ \int_{\mathbb{Z}} g(t, \hat{X}^{(j)}_t, \mathcal{L}_{\bar X_t},z)(\phi_t^{(j)}-1)\nu(dz)dt.
	\end{eqnarray}	
	{We can see that for $(\psi^{(j)}, \phi^{(j)}) \in (S_{1}^{m}\times S_{2}^{m}) $,  there exists a unique solution  $\hat{X}^{(j)} \in C([0,T], \mathbb{R}^{n})$ to the above  equation (\ref{3-1}). 
		Then it is easy to check  that the averaged coefficients also satisfy the linear growth condition and Lipschitz condition.
		So we can see that  $\big\{ {\hat{X}^{(j)}}\big\}_{ j\ge 1}$ is a family of equicontinuous and uniformly bounded functions
		in $C([0, T],\mathbb{R}^d)$.  Therefore, according to  the Arzel\`a-Ascoli theorem, the family $\{\hat X^{(j)}\}_{ j\ge 1}$ is pre-compact in $ C([0,T],\mathbb{R}^d)$. There exists a subsequence weakly converges to some limit, then we let $\hat X$ be any limit point. Then, there is  a subsequence of $\{\hat X^{(j)}\}_{ j\ge 1}$ ( which will be denoted by the same symbol) weakly converges to $\hat X$ in $ C([0,T],\mathbb{R}^d)$.	By taking same manner in \cite[Propsition 5.8]{W2022Large}, it is not too difficult to see  that the limit point $\tilde X$ satisfies the  ODEs \eqref{3}.}

	\textbf{Step 2}.  According to the variational representation for McKean-Vlasov system \cite[Theorem 3.8]{W2022Large}, we could give  the following controlled system related to (\ref{1}) as following
	\begin{eqnarray}\label{2}
	\left
	\{
	\begin{array}{ll}
	d\hat{X}^{\varepsilon ,\delta}_t =& b_1(\hat{X}^{\varepsilon ,\delta}_t,\mathcal{L}_{X^{\varepsilon ,\delta}_t}, \hat{Y}^{\varepsilon ,\delta}_t)dt + \sigma_{1}(\hat{X}^{\varepsilon ,\delta}_t,\mathcal{L}_{X^{\varepsilon ,\delta}_t})\psi_t^{\varepsilon ,\delta}dt+ \sqrt \varepsilon  \sigma_{1}(\hat{X}^{\varepsilon ,\delta}_t,\mathcal{L}_{X^{\varepsilon ,\delta}_t})dW_{t}\\
	&+ \varepsilon \int_{\mathbf{X}} g(t,\hat{X}^{\varepsilon ,\delta}_t,\mathcal{L}_{X^{\varepsilon ,\delta}_t}, z)[N^{\frac{\phi_t^{\varepsilon ,\delta}}{\varepsilon }}(dzdt)-\nu(dz)\times \frac{1}{\varepsilon}dt],\\
	 d\hat{Y}^{\varepsilon ,\delta}_t =& \frac{1}{\delta} b_2( \hat{X}^{\varepsilon ,\delta}_t,\mathcal{L}_{X^{\varepsilon ,\delta}_t}, \hat{Y}^{\varepsilon ,\delta}_t)dt +{\frac{ 1}{\sqrt{\varepsilon \delta}}}\sigma_2(\hat{X}^{\varepsilon ,\delta}_t,\mathcal{L}_{X^{\varepsilon ,\delta}_t}, \hat{Y}^{\varepsilon ,\delta}_t)\psi_t^{\varepsilon ,\delta}dt+ {\frac{ 1}{\sqrt{ \delta}}}\sigma_2(\hat{X}^{\varepsilon ,\delta}_t,\mathcal{L}_{X^{\varepsilon ,\delta}_t}, \hat{Y}^{\varepsilon ,\delta}_t)dW_{t},
	\end{array}
	\right.
	\end{eqnarray}
	where $(\psi^{\varepsilon ,\delta}, \phi^{\varepsilon ,\delta})\in \mathcal{U}^{m}$ is so-called a pair of control, according to  \cite[Theorem 3.8]{W2022Large}, it is not too difficult to see that there exists a unique solution $(\hat{X}^{\varepsilon ,\delta}, \hat{Y}^{\varepsilon ,\delta})$ to the controlled system (\ref{2}) in $\mathbf{D}\times \mathbf{D}$.

Assume that $(\psi^{\varepsilon ,\delta}, \phi^{\varepsilon ,\delta})\in \mathcal{U}^{m}$ such that $(\psi^{\varepsilon ,\delta}, \phi^{\varepsilon ,\delta})$ converges weakly to $(u, v)$ as $\varepsilon \to 0$. Then we rewrite  the slow variables in controlled  system (\ref{2}) as following,
\[
\hat{X}_{t}^{\varepsilon ,\delta}:=\mathcal{G}^{\varepsilon ,\delta}(\sqrt \varepsilon W_{t}+\int_{0}^{t}\psi_s^{\varepsilon ,\delta}ds, \varepsilon N^{\frac{\phi_t^{\varepsilon ,\delta}}{\varepsilon }}).
\]
In the following proof  in  Step 2, we will prove  that as $\varepsilon\rightarrow 0$, $\hat{X}^{\varepsilon ,\delta}$ weakly converges  to $\hat{X}$ (converges in distribution), that is
\begin{eqnarray} \label{step3}
\mathcal{G}^{\varepsilon ,\delta}(\sqrt \varepsilon W_{t}+\int_{0}^{t}\psi_s^{\varepsilon ,\delta}ds, \varepsilon N^{\frac{\phi_t^{\varepsilon ,\delta}}{\varepsilon }})\xrightarrow{\textrm{weakly}}\mathcal{G}^{0}(u, v).
\end{eqnarray}	
				
				Before showing \eqref{step3} holds, it suffices to make the following preliminary observation, there exists some  constant $C>0$ which is independent of $\varepsilon,\delta$ such that 
      \begin{eqnarray}\label{lemma3.15}
		\mathbb{E}\big[ {\mathop {\sup }\limits_{0 \le t \le T} | \hat{X}_{t}^{\varepsilon ,\delta} |^2} \big]
		<C ,\qquad
		\int_{0}^{T}\mathbb{E}\big[ { | \hat{Y}_{t}^{\varepsilon ,\delta} |^2} \big]dt<C.
		\end{eqnarray}
				
		Using the It\^o's formula directly, we get
		\begin{eqnarray}\label{lemma3.10}
		\begin{aligned}
		\mathbb{E}[{| \hat{Y}_{t}^{\varepsilon ,\delta} |^2}] &= 	\mathbb{E}[{| Y_0 |^2}] + \frac{2}{\delta }	\mathbb{E}\big[\int_0^t {\langle \hat{Y}_{s}^{\varepsilon ,\delta},b_2( \hat{X}_{s}^{\varepsilon ,\delta},\mathcal{L}_{X^{\varepsilon ,\delta}_s},\hat{Y}_{s}^{\varepsilon ,\delta} )\rangle ds}\big] \\
		&\quad
		+ \frac{2}{{\sqrt {\delta \varepsilon } }}	\mathbb{E}\big[\int_0^t {\langle \hat{Y}_{s}^{\varepsilon ,\delta},\sigma_{2} ( {\hat{X}_{s}^{\varepsilon ,\delta},\mathcal{L}_{X^{\varepsilon ,\delta}_s},\hat{Y}_{s}^{\varepsilon ,\delta}} ){\psi_s^{\varepsilon ,\delta} }\rangle ds}\big]   \\
		&\quad+\frac{2}{\sqrt{\delta} }	\mathbb{E}\big[\int_0^t {\langle \hat{Y}_{s}^{\varepsilon ,\delta},\sigma_{2}( {\hat{Y}_{s}^{\varepsilon ,\delta},\mathcal{L}_{X^{\varepsilon ,\delta}_s},\hat{Y}_{s}^{\varepsilon ,\delta}} ) dW_s\rangle}\big] \\
		&\quad+ \frac{1}{\delta }	\mathbb{E}\big[\int_0^t {|{\sigma_{2}}( {\hat{X}_{s}^{\varepsilon ,\delta},\mathcal{L}_{X^{\varepsilon ,\delta}_s},\hat{Y}_{s}^{\varepsilon ,\delta}} )|^2ds}\big] .
		\end{aligned}
		\end{eqnarray}

		{The fourth term is a true martingale. In particular, we have $\mathbb{E}[\int_0^t {\langle \hat {Y}_{s}^{\varepsilon ,\delta},\sigma_{2}( {\hat {X}_{s}^{\varepsilon ,\delta},\mathcal{L}_{X^{\varepsilon ,\delta}_s},\hat {Y}_{s}^{\varepsilon ,\delta}} ) dW_s\rangle}]=0$.
		So, we can get that
		\begin{eqnarray}\label{3-41}
		\begin{aligned}
		\frac{d\mathbb{E}[{| \hat {Y}_{t}^{\varepsilon ,\delta} |^2}]}{dt} &= \frac{2}{\delta }\mathbb{E} \big[{\langle \hat {Y}_{t}^{\varepsilon ,\delta},b_2( \hat {X}_{t}^{\varepsilon ,\delta},\mathcal{L}_{X^{\varepsilon ,\delta}_t},\tilde {Y}_{t}^{\varepsilon ,\delta} )\rangle } \big]
		+ \frac{2}{{\sqrt {\delta \varepsilon } }}\mathbb{E} \big[{\langle \hat {Y}_{t}^{\varepsilon ,\delta},\sigma_{2} ( {\hat {X}_{t}^{\varepsilon ,\delta},\mathcal{L}_{X^{\varepsilon ,\delta}_t},\mathcal{L}_{X^{\varepsilon ,\delta}_t},\hat {Y}_{t}^{\varepsilon ,\delta}} )\psi_t^{\varepsilon ,\delta}\rangle } \big]  \\
		&\quad+ \frac{1}{\delta }\mathbb{E}\big[ |{\sigma_{2} }( {\hat {X}_{t}^{\varepsilon ,\delta},\mathcal{L}_{X^{\varepsilon ,\delta}_t},\hat {Y}_{t}^{\varepsilon ,\delta}} )|^2\big].
		\end{aligned}
		\end{eqnarray} 
		With Assumption (\textbf{A3}), we have
		\begin{eqnarray}\label{lemma3.12}
		&\frac{2}{\delta }{\langle \hat {Y}_{t}^{\varepsilon ,\delta},b_2( \hat {X}_{t}^{\varepsilon ,\delta},\mathcal{L}_{X^{\varepsilon ,\delta}_t},\hat {Y}_{t}^{\varepsilon ,\delta} )\rangle } +\frac{1}{\delta } |{\sigma_{2} }( {\hat {X}_{t}^{\varepsilon ,\delta},\mathcal{L}_{X^{\varepsilon ,\delta}_t},\hat {Y}_{t}^{\varepsilon ,\delta}} )|^2	\cr	&\le- \frac{{ 2C_5 }}{\delta }{{| \hat{Y}_{t}^{\varepsilon ,\delta} |^2}}  + \frac{{C_6}}{\delta }(1+{| \hat {X}_{t}^{\varepsilon ,\delta} |^2}+\mathcal{L}_{X^{\varepsilon ,\delta}_t}(|\cdot|^2)) .	
		\end{eqnarray}
		By  Assumption (\textbf{A2}) and the fact that  $(\psi^{\varepsilon ,\delta}, \phi^{\varepsilon ,\delta})\in \mathcal{U}^{m}$,  we have
		\begin{eqnarray}\label{lemma3.13}
		\begin{aligned}
		\frac{2}{{\sqrt {\delta \varepsilon } }}{\langle \hat {Y}_{t}^{\varepsilon ,\delta},\sigma_{2} ( {\hat {X}_{t}^{\varepsilon ,\delta},\mathcal{L}_{X^{\varepsilon ,\delta}_t},\tilde {Y}_{t}^{\varepsilon ,\delta}} )\psi_t^{\varepsilon ,\delta}\rangle } 			&\le   \frac{{{1 }}}{{\sqrt {\delta \varepsilon } }} {  {| \hat{Y}_{t}^{\varepsilon ,\delta} |^2} }+\frac{{{C_3 }}}{{\sqrt {\delta \varepsilon } }} {\big( {1 + {| \hat{X}_{t}^{\varepsilon ,\delta} |^2}+\mathcal{L}_{X^{\varepsilon ,\delta}_t}(|\cdot|^2)} \big){(\psi_t^{\varepsilon ,\delta} )^2}} \\ .
		\end{aligned}
		\end{eqnarray}
		Thus, as a consequence of  (\ref{lemma3.10})-- (\ref{lemma3.13}), it deduces  that
		\begin{eqnarray}\label{lemma3.17}
		\begin{aligned}
		\frac{d\mathbb{E}[{| \hat{Y}_{t}^{\varepsilon ,\delta} |^2}]}{dt} &\le  {\frac{{ - 2C_5 }}{\delta } } 	\frac{d\mathbb{E}[{| \hat{Y}_{t}^{\varepsilon ,\delta} |^2}]}{dt} +\frac{{{1 }}}{{\sqrt {\delta \varepsilon } }}	\frac{d\mathbb{E}[{| \hat{Y}_{t}^{\varepsilon ,\delta} |^2}]}{dt} + \frac{{{C_3 }}}{{\sqrt {\delta \varepsilon } }} 	\mathbb{E}[{{| \hat{X}_{t}^{\varepsilon ,\delta} |^2{| {\psi_t^{\varepsilon ,\delta}} |^2}}}]  \\
		&\quad +   \frac{{{C_3 }\mathbb{E}[{|\psi_t^{\varepsilon ,\delta} |^2}]}}{{\sqrt {\delta \varepsilon } }}  + \frac{{2 C_6}}{\delta }+\frac{{{C_3 }{|\psi_t^{\varepsilon ,\delta} |^2}}}{{\sqrt {\delta \varepsilon } }}\sup_{t\in[0,T]}\mathbb{E}[|X_t^{\varepsilon ,\delta}|^2].
		\end{aligned}
		\end{eqnarray}
		So, by comparison theorem, we have that
		\begin{eqnarray*}\label{3-231}
			\begin{aligned}
				\mathbb{E}[{| \hat {Y}_{t}^{\varepsilon ,\delta} |^2}] &\le  |y_0|^2 e^{-\frac{2C_5}{\delta} t} +   \frac{{{C_3 }}}{{\sqrt {\delta \varepsilon } }}\int_{0}^{t}  e^{-\frac{2C_5}{\delta} (t-s)}{\mathbb{E}[{| \hat {X}_{s}^{\varepsilon ,\delta} |^2{| \psi_s^{\varepsilon ,\delta} |^2}}]} ds +   \frac{{{C_3 }}}{\sqrt {\delta \varepsilon } }\int_{0}^{t}  e^{-\frac{2C_5}{\delta} (t-s)}\mathbb{E}[{|\psi_s^{\varepsilon ,\delta} |^2}] ds  \\
				&\qquad+ \frac{{ 2C_6}}{\delta }\int_{0}^{t}  e^{-\frac{2C_5}{\delta} (t-s)}ds+\frac{{{C_3 }}}{{\sqrt {\delta \varepsilon } }} \sup_{t\in[0,T]}\mathbb{E}[|X_t^{\varepsilon ,\delta}|^2]\int_{0}^{t}e^{-\frac{2C_5}{\delta} (t-s)}{|\psi_s^{\varepsilon ,\delta} |^2}ds.
			\end{aligned}
		\end{eqnarray*}
		Then, by using Fubini theorem and \lemref{lem1}, it deduces that
		\begin{eqnarray*}\label{3-232}
			\begin{aligned}
				\int_{0}^{T}\mathbb{E}[{| \hat {Y}_{t}^{\varepsilon ,\delta} |^2}]dt 
				&\le  |y_0|^2 \int_{0}^{T}e^{-\frac{2C_5}{\delta} t}dt +   \frac{{{C_3 }}}{{\sqrt {\delta \varepsilon } }}\int_{0}^{T}\int_{0}^{t}  e^{-\frac{2C_5}{\delta} (t-s)}{\mathbb{E}[{| \hat {X}_{s}^{\varepsilon ,\delta} |^2{| \psi_s^{\varepsilon ,\delta} |^2}}]} dsdt\\
				&\quad +   \frac{{{C_3 }}}{\sqrt {\delta \varepsilon } }\int_{0}^{T}\int_{0}^{t}  e^{-\frac{2C_5}{\delta} (t-s)} {|\psi_s^{\varepsilon ,\delta} |^2}dsdt+ \frac{{ 2C_6}}{\delta }\int_{0}^{T}\int_{0}^{t}  e^{-\frac{2C_5}{\delta} (t-s)}dsdt  \\
				&\quad+\frac{{{C_3 }}}{{\sqrt {\delta \varepsilon } }} \sup_{t\in[0,T]}\mathbb{E}[|X_t^{\varepsilon ,\delta}|^2]\int_{0}^{T}\int_{0}^{t}e^{-\frac{2C_5}{\delta} (t-s)}{|\psi_s^{\varepsilon ,\delta} |^2}dsdt\\
				&\le |y_0|^2 e^{-\frac{2C_5}{\delta} T} +   \frac{\delta C_3}{C_5\sqrt {\delta \varepsilon } }\mathbb{E}\big[\sup_{0 \leq t \leq T}| \hat {X}_{t}^{\varepsilon ,\delta} |^2\int_{0}^{T}  e^{-\frac{2C_5}{\delta} (T-s)} | \psi_s^{\varepsilon ,\delta} |^2ds\big]   + C.
			\end{aligned}
		\end{eqnarray*}
		With aid of the fact that $(\psi^{\varepsilon ,\delta}, \phi^{\varepsilon ,\delta})\in \mathcal{U}^m$, it deduces 
		\begin{eqnarray*}\label{lemma3.171}
			\int_{0}^{T}\mathbb{E}[{| \hat {Y}_{t}^{\varepsilon ,\delta} |^2}]dt &\le& {C}\mathbb{E}[\mathop {\sup }\limits_{t \in \left[ {0,T} \right]} {| \hat {X}_{t}^{\varepsilon ,\delta} |^2}]+C.
		\end{eqnarray*}}
	Likewise, by using the It\^o's formula, we get
		\begin{eqnarray}\label{lemma3.211}
		{| \hat{X}_{t}^{\varepsilon ,\delta} |^2}={| {x_0} |^2}  +{\mathcal{I}_1} + {\mathcal{I}_2} + {\mathcal{I}_3} + {\mathcal{I}_4} + {\mathcal{I}_{51}}+{\mathcal{I}_{52}} + {\mathcal{I}_6} + {\mathcal{I}_7} ,
		\end{eqnarray}
		where 
		\begin{eqnarray*}
		\begin{aligned}
		 &{\mathcal{I}_1}= 2\int_0^t {\langle \hat{X}_{s}^{\varepsilon ,\delta},b_1( {\hat{X}_{s}^{\varepsilon ,\delta},\mathcal{L}_{X^{\varepsilon ,\delta}_s},\hat{Y}_{s}^{\varepsilon ,\delta}} )\rangle ds},\quad
		 	{\mathcal{I}_2}= 2\sqrt \varepsilon  \int_0^t {\langle \hat{Y}_{s}^{\varepsilon ,\delta},\sigma_{1}( \hat{X}_{s}^{\varepsilon ,\delta},\mathcal{L}_{X^{\varepsilon ,\delta}_s} ) dW_s\rangle} ,\\
		&{\mathcal{I}_3}= 2\int_0^t {\langle \hat{X}_{s}^{\varepsilon ,\delta},\sigma_{1}( \hat{X}_{s}^{\varepsilon ,\delta},\mathcal{L}_{X^{\varepsilon ,\delta}_s} )\psi_s^{\varepsilon ,\delta}\rangle ds},\quad
		{\mathcal{I}_4}= \varepsilon \int_0^t {|{\sigma_{1}}( \hat{X}_{s}^{\varepsilon ,\delta} ,\mathcal{L}_{X^{\varepsilon ,\delta}_s})|^2ds},\\
		&{\mathcal{I}_{51}}= \int_0^t {\int_{\mathbf{X}} { {{\varepsilon^2 g^2(s, \hat{X}_{s}^{\varepsilon ,\delta} ,\mathcal{L}_{X^{\varepsilon ,\delta}_s},z) } } {\tilde N}^{{\phi_s^{\varepsilon ,\delta} } \mathord{\left/{\vphantom {{\phi_s^{\varepsilon ,\delta} } \varepsilon }} \right.\kern-\nulldelimiterspace} \varepsilon }\left( {dzds} \right)} },\\
		&{\mathcal{I}_{52}} = \int_0^t {\int_{\mathbf{X}} [ {2\varepsilon \hat{X}_{s}^{\varepsilon ,\delta}g(s, \hat{X}_{s}^{\varepsilon ,\delta} ,\mathcal{L}_{X^{\varepsilon ,\delta}_s}, z)} ]{\tilde N}^{{\phi_s^{\varepsilon ,\delta} } \mathord{\left/{\vphantom {{\phi_s^{\varepsilon ,\delta} } \varepsilon }} \right.\kern-\nulldelimiterspace} \varepsilon }( {dzds} ) },\\
		&{\mathcal{I}_6}= \int_0^t {\int_{\mathbf{X}} {\varepsilon |{g}(s, \hat{X}_{s}^{\varepsilon ,\delta} ,\mathcal{L}_{X^{\varepsilon ,\delta}_s},z)|^2\phi_s^{\varepsilon ,\delta} \nu(dz)ds} } ,\\
		&{\mathcal{I}_7}=\int_0^t {\int_{\mathbf{X}} {2\langle \hat{X}_{s}^{\varepsilon ,\delta},g(s, \hat{X}_{s}^{\varepsilon ,\delta},\mathcal{L}_{X^{\varepsilon ,\delta}_s},z )(\phi_s^{\varepsilon ,\delta}-1)\rangle \nu(dz)ds} }.
		\end{aligned}
		\end{eqnarray*}
According to Assumptions (\textbf{A1}), (\textbf{A4})  and some straightforward computation, we have  following estimates,
		\begin{eqnarray}\label{lemma3.3}
		\begin{aligned}
			\left| {{\mathcal{I}_1}} \right| &\le {C_2}\int_0^t (1+{| \hat {X}_{s}^{\varepsilon ,\delta} |^2}+\mathcal{L}_{X^{\varepsilon ,\delta}_s}(|\cdot|^2))ds+\int_0^t {| \hat {X}_{s}^{\varepsilon ,\delta} |^2}ds   + {C_2}\int_0^t {| \hat{Y}_{s}^{\varepsilon ,\delta} |^2ds},\\	
			\left| {{\mathcal{I}_3}} \right| &\le  \int_0^t {| \hat{X}_{s}^{\varepsilon ,\delta} |^2 | {\psi_s^{\varepsilon ,\delta} } |^2ds}  + {C_2}  {\int_0^t (1+{| \hat {X}_{s}^{\varepsilon ,\delta} |^2}+\mathcal{L}_{X^{\varepsilon ,\delta}_s}(|\cdot|^2))ds} \\
			&\le M_\psi\int_0^t { | \hat{X}_{s}^{\varepsilon ,\delta} |^2ds}  + {C_2}\int_0^t { | \hat{X}_{s}^{\varepsilon ,\delta} |^2ds} +C_2\sup_{t\in[0,T]}\mathbb{E}[|X_t^{\varepsilon ,\delta}|^2] + {C_2}T,\\
			\left| {{\mathcal{I}_4}} \right| &\le \varepsilon  {\int_0^t [{| \hat {X}_{s}^{\varepsilon ,\delta} |^2}+\mathcal{L}_{X^{\varepsilon ,\delta}_s}(|\cdot|^2)]ds}  + \varepsilon {C_2}T,\\
			\left| {{\mathcal{I}_6}} \right| &\le \varepsilon {C_2} \Big\{ {\int_0^t {\int_{\mathbf{X}} {(1+{| \hat {X}_{s}^{\varepsilon ,\delta} |^2}+\mathcal{L}_{X^{\varepsilon ,\delta}_s}(|\cdot|^2)){\|g(s, z)\|}\phi^{\varepsilon ,\delta} \nu\left( {dz} \right)ds} } } \Big\}\\
			&\le \varepsilon {C_2}M_\phi + \varepsilon {C_2}M_\phi \mathop {\sup }\limits_{t \in \left[ {0,T} \right]} | \hat{X}_{t}^{\varepsilon ,\delta} |^2+ \varepsilon {C_2}M_\phi T \mathop {\sup }\limits_{t \in \left[ {0,T} \right]} \mathbb{E}[| {X}_{t}^{\varepsilon ,\delta} |^2],\\
			\left| {{\mathcal{I}_7}} \right| &\le   {\int_0^t {\int_{\mathbf{X}} | \hat{X}_{s}^{\varepsilon ,\delta} |^2(\phi_s^{\varepsilon ,\delta}-1)\nu( dz ) ds} } \\
		&\quad\quad+ {C_2} {\int_0^t {\int_{\mathbf{X}} (1+{| \hat {X}_{s}^{\varepsilon ,\delta} |^2}+\mathcal{L}_{X^{\varepsilon ,\delta}_s}(|\cdot|^2)){\|g(s, z)\|}(\phi_s^{\varepsilon ,\delta}-1)\nu\left( {dz} \right)ds } },
		\end{aligned}
		\end{eqnarray}
				where 		 $M_\psi=\mathop {\sup }\limits_{\psi \in S_{2}^{m}}\int_0^T { ( \psi_s^{\varepsilon ,\delta} )^2ds}<\infty$,
and 
\[
	M_\phi=\max\{\mathop {\sup }\limits_{\phi \in  S_1^{m}}{\int_0^t {\int_{\mathbf{X}} | {\phi_s^{\varepsilon ,\delta}  - 1} |\nu( dz )ds } },\mathop {\sup }\limits_{\phi \in  S_1^{m}}{\int_0^t {\int_{\mathbf{X}} {{\|g(s, z)\|}| {\phi_s^{\varepsilon ,\delta}  - 1} |\nu( {dz} )ds} } }\}<\infty,
\]
which is deduced from \cite[Lemma 3.4]{2013Large}.
Then with the Gronwall's inequality and \lemref{lem1}, it deduces
		\begin{eqnarray} \label{lemma3.212}
			\mathop {\sup }\limits_{t \in \left[ {0,T} \right]}{| \hat{X}_{t}^{\varepsilon ,\delta} |^2} &\le& c_{15}\big({| x_0|^2}+\mathop {\sup }\limits_{t \in \left[ {0,T} \right]} \left| {{\mathcal{I}_2}} \right|+\mathop {\sup }\limits_{t \in \left[ {0,T} \right]} | {{\mathcal{I}_{51}}} |+\mathop {\sup }\limits_{t \in \left[ {0,T} \right]} | {{\mathcal{I}_{52}}} |+ {C}T + {C_2}\int_0^T {| \hat{Y}_{s}^{\varepsilon ,\delta} |^2ds}\big).
		\end{eqnarray}				
By the H\"older inequality and Assumption (\textbf{A1}), it follows 
		\begin{eqnarray}\label{lemma3.4}
		\begin{aligned}
		\mathbb{E}[\mathop {\sup }\limits_{t \in \left[ {0,T} \right]} \left| {{\mathcal{I}_2}} \right|] &\le 4\mathbb{E}{\Big[ {4\varepsilon \int_0^T {| \hat{X}_{s}^{\varepsilon ,\delta} |^2|\sigma_{1}( \hat{X}_{s}^{\varepsilon ,\delta},\mathcal{L}_{X^{\varepsilon ,\delta}_s} )|^2ds} } \Big]^{{1 \mathord{\left/{\vphantom {1 2}} \right.\kern-\nulldelimiterspace} 2}}}\\
		&\le 8{C_2}\mathbb{E}{\Big[ {\varepsilon \mathop {\sup }\limits_{t \in \left[ {0,T} \right]} {{| \hat{X}_{t}^{\varepsilon ,\delta} |}^2}\int_0^T {(1+{| \hat {X}_{s}^{\varepsilon ,\delta} |^2}+\mathcal{L}_{X^{\varepsilon ,\delta}_s}(|\cdot|^2))ds} } \Big]^{{1 \mathord{\left/{\vphantom {1 2}} \right.\kern-\nulldelimiterspace} 2}}}\\
		&\le \frac{1}{8c_{15}}\mathbb{E}[\mathop {\sup }\limits_{t \in \left[ {0,T} \right]} {| \hat{X}_{t}^{\varepsilon ,\delta} |^2}] + 128 c_{15}\varepsilon\mathbb{E}\Big[ \int_0^T (1+{| \hat {X}_{s}^{\varepsilon ,\delta} |^2}+\mathcal{L}_{X^{\varepsilon ,\delta}_s}(|\cdot|^2))ds  \Big]\\
		&\le \frac{1}{8c_{15}}\mathbb{E}[\mathop {\sup }\limits_{t \in \left[ {0,T} \right]} {| \hat{X}_{t}^{\varepsilon ,\delta} |^2}] + 128c_{15}\varepsilon T\\
		&\quad+ 128c_{15}\varepsilon T \mathbb{E}[\mathop {\sup }\limits_{t \in \left[ {0,T} \right]} {| \hat{X}_{t}^{\varepsilon ,\delta} |^2}]+ 128c_{15}\varepsilon T \mathbb{E}[\mathop {\sup }\limits_{t \in \left[ {0,T} \right]} {| {X}_{t}^{\varepsilon ,\delta} |^2}].
		\end{aligned}
		\end{eqnarray}
	What is left is to estimate remaining terms ${{\mathcal{I}_{51}}}$ and ${{\mathcal{I}_{52}}}$.	By Assumption (\textbf{A4}) and the Burkholder-Davis-Gundy inequality, it has
		\begin{eqnarray}\label{lemma3.7}
		\begin{aligned}
		\mathbb{E}[\mathop {\sup }\limits_{t \in \left[ {0,T} \right]} \left| {{\mathcal{I}_{51}}} \right|] &\le 2{\varepsilon }\mathbb{E}\Big[ {\int_0^T {\int_{\mathbf{X}} {{g^2}(s, \hat{X}_{s}^{\varepsilon ,\delta},\mathcal{L}_{X^{\varepsilon ,\delta}_s},z )\phi_s^{\varepsilon ,\delta}\nu\left( {dz} \right)ds} } } \Big] \\
		&\le 4\varepsilon {C_2}M_\phi\mathbb{E}[{\mathop {\sup }\limits_{t \in \left[ {0,T} \right]} | \hat{X}_{t}^{\varepsilon ,\delta} |^2}]+  4\varepsilon {C_2}M_\phi\mathbb{E}[{\mathop {\sup }\limits_{t \in \left[ {0,T} \right]} | {X}_{t}^{\varepsilon ,\delta} |^2}]+ 4\varepsilon {C_2}M_\phi.\\
		\mathbb{E}[\mathop {\sup }\limits_{t \in \left[ {0,T} \right]} \left| {{\mathcal{I}_{52}}} \right|] &\le 4\mathbb{E}\left[ {\mathcal{I}_{_{52}}^2} \right]_T^{{1 \mathord{\left/{\vphantom {1 2}} \right.\kern-\nulldelimiterspace} 2}}\\
		&\le 4\mathbb{E}{\Big[{4{\varepsilon ^2}\int_0^T {\int_{\mathbf{X}} {{( \hat{X}_{s}^{\varepsilon ,\delta} )^2}{g^2}(s, \hat{X}_{s}^{\varepsilon ,\delta},\mathcal{L}_{X^{\varepsilon ,\delta}_s},z ){{ N}^{{{\phi_s^{\varepsilon ,\delta}} \mathord{\left/	{\vphantom {{\phi^{\varepsilon ,\delta} } \varepsilon }} \right.\kern-\nulldelimiterspace} \varepsilon }}}\left( {dzds} \right)} } } \Big]^{{1 \mathord{\left/{\vphantom {1 2}} \right.\kern-\nulldelimiterspace} 2}}}\\
		&\le 8\mathbb{E}{\Big[ {\varepsilon^2 \mathop {\sup }\limits_{t \in \left[ {0,T} \right]} {\big| \hat{X}_{t}^{\varepsilon ,\delta} \big|^2}\int_0^T {\int_{\mathbf{X}} {{g^2}(s, \hat{X}_{s}^{\varepsilon ,\delta},\mathcal{L}_{X^{\varepsilon ,\delta}_s},z ){{ N}^{{{\phi_s^{\varepsilon ,\delta}} \mathord{\left/	{\vphantom {{\phi^{\varepsilon ,\delta} \left( t \right)} \varepsilon }} \right.\kern-\nulldelimiterspace} \varepsilon }}}(dzds)} } } \Big]^{{1 \mathord{\left/{\vphantom {1 2}} \right.\kern-\nulldelimiterspace} 2}}}\\
		&\le \frac{1}{8c_{15}}\mathbb{E}[\mathop {\sup }\limits_{t \in \left[ {0,T} \right]} {| \hat{X}_{t}^{\varepsilon ,\delta} |^2}] + 128c_{16} M_\phi\varepsilon\mathbb{E}[{\mathop {\sup }\limits_{t \in \left[ {0,T} \right]} | \hat{X}_{t}^{\varepsilon ,\delta} |^2}] \\
		&\quad+ 128c_{16} M_\phi\varepsilon\mathbb{E}[{\mathop {\sup }\limits_{t \in \left[ {0,T} \right]} | {X}_{t}^{\varepsilon ,\delta} |^2}] + 128c_{16}M_\phi\varepsilon.
		\end{aligned}
		\end{eqnarray}
		With the help of  the Gronwall's inequality,  estimates (\ref{lemma3.212})--(\ref{lemma3.7}), we  obtain that
		\begin{eqnarray}\label{lemma3.102}
			\mathbb{E}[\mathop {\sup }\limits_{t \in \left[ {0,T} \right]} {| \hat{X}_{t}^{\varepsilon ,\delta} |^2}] \le c_{17}.
		\end{eqnarray}		
Thus, estimates (\ref{lemma3.15}) can be obtained.

Next, it remains to show (\ref{step3}) as $\varepsilon\rightarrow 0$. Firstly, we  construct the following auxiliary processes. Set $t(\Delta):={\left[ \frac{t}{ \Delta}\right] }\Delta$, then define
	\begin{eqnarray*}
	\left
	\{
	\begin{array}{ll}
	d\tilde{X}_{t}^{\varepsilon ,\delta} &=  b_1( {{ \hat{X}_{t(\Delta)}^{\varepsilon ,\delta}},\mathcal{L}_{X^{\varepsilon ,\delta}_{t(\Delta)}},\tilde y_t^{\varepsilon ,\delta }} )dt+ \sigma_{1}( \tilde{X}_{t}^{\varepsilon ,\delta} ,\mathcal{L}_{X^{\varepsilon ,\delta}_{t}})\psi_t^{\varepsilon ,\delta} dt + \int_{\mathbf{X}} {g(t, {\tilde{X}_{t}^{\varepsilon ,\delta}},\mathcal{L}_{X^{\varepsilon ,\delta}_{t}},z )( {\phi_t^{\varepsilon ,\delta} - 1} )\nu(dz)dt}, \\
	d\tilde Y_t^{\varepsilon ,\delta } &=   \frac{1}{\delta }b_2\big( {{ \hat{X}_{t(\Delta)}^{\varepsilon ,\delta}},\mathcal{L}_{X^{\varepsilon ,\delta}_{t(\Delta)}},\tilde y_t^{\varepsilon ,\delta }} \big)dt + \frac{1}{{\sqrt \delta  }}\sigma_{2} \big( {{ \hat{X}_{t(\Delta)}^{\varepsilon ,\delta}},\mathcal{L}_{X^{\varepsilon ,\delta}_{t(\Delta)}},\tilde Y_t^{\varepsilon ,\delta }} \big)dW_t ,
	\end{array}
	\right.
	\end{eqnarray*}
{by taking the same manner in \eqref{lemma3.10}--\eqref{lemma3.102}}, it deduces that 
		\begin{eqnarray}\label{3-2}
			\mathbb{E}\big[ {\mathop {\sup }\limits_{0 \le t \le T} | \tilde{X}_{t}^{\varepsilon ,\delta} |^2} \big]< \infty,\qquad
			\mathbb{E}\big[ {| {\tilde Y_t^{\varepsilon ,\delta }} |^2} \big]< \infty.
		\end{eqnarray}
Then, we  construct the  stopping times as follows, for any $R, \varepsilon>0$, set
		{$\tau_R^\varepsilon:=inf\{t\in[0,T]:| \hat{X}_{t}^{\varepsilon ,\delta} |>R\}.$}
For $t,t-h \in [0,  T\wedge\tau_R^\varepsilon]$, it will show that
		\begin{eqnarray}\label{lemma3.18}
		\mathbb{E}[\int_{0}^{T\wedge\tau_R^\varepsilon}|\hat{X}_{t}^{\varepsilon ,\delta} - { \hat{X}_{t-h}^{\varepsilon ,\delta}} |^2dt]\le Ch.
		\end{eqnarray}
		Indeed, it implies from It\^o's formula directly that
		\begin{eqnarray}\label{lemma3.21}
		{\int_{0}^{T\wedge\tau_R^\varepsilon}|\hat{X}_{t}^{\varepsilon ,\delta} - { \hat{X}_{t-h}^{\varepsilon ,\delta}}|^2}dt
		= {{\mathcal{H}}_1} + {{\mathcal{H}}_2} + {{\mathcal{H}}_3} + {{\mathcal{H}}_4} + {{\mathcal{H}}_5} + {{\mathcal{H}}_6} + {{\mathcal{H}}_7},
		\end{eqnarray}
		where
		\begin{eqnarray*}
		\begin{aligned}
       {{\mathcal{H}}_1}
		&=  2\int_{0}^{T\wedge\tau_R^\varepsilon}\int_{t-h}^t {\langle ( {\hat{X}_{s}^{\varepsilon ,\delta} - \hat{X}_{t-h}^{\varepsilon ,\delta}} ),b_1( {\hat{X}_{s}^{\varepsilon ,\delta},\mathcal{L}_{X^{\varepsilon ,\delta}_s},\hat{Y}_{s}^{\varepsilon ,\delta}} )\rangle ds}dt, \\
		{{\mathcal{H}}_2}
		&=  2\sqrt \varepsilon  \int_{0}^{T\wedge\tau_R^\varepsilon}\int_{t-h}^t {\big\langle ( {\hat{X}_{s}^{\varepsilon ,\delta} - \hat{X}_{t-h}^{\varepsilon ,\delta}} ),\sigma_{1}( \hat{X}_{s}^{\varepsilon ,\delta} ,\mathcal{L}_{X^{\varepsilon ,\delta}_s}) dW_s\big\rangle}dt,\\
		{{\mathcal{H}}_3}
		&=  2\int_{0}^{T\wedge\tau_R^\varepsilon}\int_{t-h}^t {\langle ( {\hat{X}_{s}^{\varepsilon ,\delta} - \hat{X}_{t-h}^{\varepsilon ,\delta}} ),\sigma_{1}( \hat{X}_{s}^{\varepsilon ,\delta},\mathcal{L}_{X^{\varepsilon ,\delta}_s} )\psi_s^{\varepsilon ,\delta} \rangle ds}dt, \\
		{{\mathcal{H}}_4}
		&=  \varepsilon \int_{0}^{T\wedge\tau_R^\varepsilon}\int_{t-h}^t {|{\sigma_{1}}( \hat{X}_{s}^{\varepsilon ,\delta} ,\mathcal{L}_{X^{\varepsilon ,\delta}_s})|^2ds}dt,\\
		{{\mathcal{H}}_5}
		&=  \int_{0}^{T\wedge\tau_R^\varepsilon}\int_{t-h}^t {\int_{\mathbf{X}} {\big[ {{( {( {\hat{X}_{s}^{\varepsilon ,\delta} - \hat{X}_{t-h}^{\varepsilon ,\delta}} ){\rm{ + }}\varepsilon g(s, \hat{X}_{s}^{\varepsilon ,\delta},\mathcal{L}_{X^{\varepsilon ,\delta}_s},z )} )^2} - {( {\hat{X}_{s}^{\varepsilon ,\delta} - \hat{X}_{t-h}^{\varepsilon ,\delta}} )^2}} \big]{{\tilde N}^{\frac{\phi_s^{\varepsilon ,\delta}}{\varepsilon}}}\left( {dzds} \right)} }dt,\\
		{{\mathcal{H}}_6}
		&=  \int_{0}^{T\wedge\tau_R^\varepsilon}\int_{t-h}^t {\int_{\mathbf{X}} {\varepsilon |g(s, \hat{X}_{s}^{\varepsilon ,\delta},\mathcal{L}_{X^{\varepsilon ,\delta}_s},z )|^2\phi_s^{\varepsilon ,\delta}\nu(dz)ds} }dt,\\
		{{\mathcal{H}}_7}
		&=  \int_{0}^{T\wedge\tau_R^\varepsilon}\int_{t-h}^t {\int_{\mathbf{X}} {2\langle ( {\hat{X}_{s}^{\varepsilon ,\delta} - \hat{X}_{t-h}^{\varepsilon ,\delta}} ),g(s, \hat {X}_{s}^{\varepsilon ,\delta},\mathcal{L}_{\bar X_s},z )\rangle(\phi_s^{\varepsilon ,\delta}-1) \nu(dz)ds} }dt.
		\end{aligned}
		\end{eqnarray*}
			According to the  Assumption {(\textbf{A2})},  H\"older inequality and Fubini theorem, we get that following estimates
		\begin{eqnarray}\label{lemma3.23}
		\begin{aligned}
		\mathbb{E}[{{\mathcal{H}}_1} ]
		&\le \mathbb{E} [\int_{0}^{T\wedge\tau_R^\varepsilon}\int_{t-h}^t {{( {\hat{X}_{s}^{\varepsilon ,\delta} - \hat{X}_{t-h}^{\varepsilon ,\delta}} )^2}ds}dt]  + {C_2}h \big(1+ R^2 +\sup_{t\in[0,T]}\mathcal{L}_{X^{\varepsilon ,\delta}_t}(|\cdot|^2) +\sup_{t\in[0,T]}\mathbb{E}[|\hat{Y}_{t}^{\varepsilon ,\delta} |^2]  \big),\\
		{{\mathcal{H}}_3} 
		&\le \int_{0}^{T\wedge\tau_R^\varepsilon}\int_{t-h}^t {{\big( {\hat{X}_{s}^{\varepsilon ,\delta} - \hat{X}_{t-h}^{\varepsilon ,\delta}} \big)^2}ds}dt  + {C_2}\int_{0}^{T\wedge\tau_R^\varepsilon}\int_{t-h}^t {\big( 1 + | \hat{X}_{s}^{\varepsilon ,\delta} |^2+\mathcal{L}_{X^{\varepsilon ,\delta}_s}(|\cdot|^2)\big) ( \psi_s^{\varepsilon ,\delta})^2ds}dt \\
		&\le \int_{0}^{T\wedge\tau_R^\varepsilon}\int_{t-h}^t {( {\hat{X}_{s}^{\varepsilon ,\delta} - \hat{X}_{t-h}^{\varepsilon ,\delta}} )^2dsdt  + {C_2}M_\psi\big( 1 +  {R^2}+\sup_{t\in[0,T]}\mathcal{L}_{X^{\varepsilon ,\delta}_t}(|\cdot|^2)} \big)h,\\
		{{\mathcal{H}}_4}   &\le 2\varepsilon {C_2}\mathbb{E}\int_{0}^{T\wedge\tau_R^\varepsilon}\int_{t-h}^t {\big( 1 + | \hat{X}_{s}^{\varepsilon ,\delta} |^2 +\mathcal{L}_{X^{\varepsilon ,\delta}_s}(|\cdot|^2)\big)ds}dt , \\
		{{\mathcal{H}}_6}  &\le \varepsilon M_\phi\big( 1 + R^2+\sup_{t\in[0,T]}\mathcal{L}_{X^{\varepsilon ,\delta}_t}(|\cdot|)^2 \big).\\
		{{\mathcal{H}}_7} &\le \int_{0}^{T\wedge\tau_R^\varepsilon}\int_{t-h}^t {\int_{\mathbf{X}} \big( {\hat{X}_{s}^{\varepsilon ,\delta} - \hat{X}_{t-h}^{\varepsilon ,\delta}} \big)^2( {\phi_s^{\varepsilon ,\delta} - 1} )\nu(dz)ds }dt\\
		&\quad+ \big( 1 +  {R^2}+\sup_{t\in[0,T]}\mathcal{L}_{X^{\varepsilon ,\delta}_t}(|\cdot|)^2 \big)\int_{0}^{T\wedge\tau_R^\varepsilon}\int_{t-h}^t {\int_{\mathbf{X}} {\|g(s, z)\|( {\phi_s^{\varepsilon ,\delta}  - 1} )\nu(dz)ds} }dt.
		\end{aligned}
		\end{eqnarray}

		From estimates (\ref{lemma3.23}), it deduces that
		\begin{eqnarray}\label{lemma3.25}
		\begin{aligned}
		\mathbb{E}[\int_{0}^{T\wedge\tau_R^\varepsilon}| \hat{X}_{t}^{\varepsilon ,\delta} - { \hat{X}_{t-h}^{\varepsilon ,\delta}}|^2 dt]	
			&\le \mathbb{E}[{{\mathcal{H}}_2} ]  +\mathbb{E}[{{\mathcal{H}}_5}] +8R^2 h+ 4R^2\int_{0}^{T\wedge\tau_R^\varepsilon}\int_{t-h}^t {\int_{\mathbf{X}} ( {\phi_s^{\varepsilon ,\delta}  - 1} )\nu(dz)ds }dt \\
			&\quad+ {C_2}h \big(  1+\mathbb{E}[\sup_{t\in[0,T]}|X_t^{\varepsilon ,\delta}|^2]+{R^2} \\
			&\quad+ \sup_{t\in[0,T]}\mathbb{E}[| \hat{Y}_{t}^{\varepsilon ,\delta} |^2] \big) + {C_2}{M_\psi }\big( {1 +  {R^2}+\mathbb{E}[\sup_{t\in[0,T]}|X_t^{\varepsilon ,\delta}|^2]} \big)h\\
			&\quad +2\varepsilon {C_2}h  (1+{R^2}+\mathbb{E}[\sup_{t\in[0,T]}|X_t^{\varepsilon ,\delta}|^2])  \\
			&\quad+\varepsilon M_\phi\big( {1 +  {{ R }^2}}+\mathbb{E}[\sup_{t\in[0,T]}|X_t^{\varepsilon ,\delta}|^2 ]\big)\\
			&\quad+ ( 1 +  {R^2}+\mathbb{E}[\sup_{t\in[0,T]}|X_t^{\varepsilon ,\delta}|^2])\int_{0}^{T\wedge\tau_R^\varepsilon}\int_{t-h}^t {\int_{\mathbf{X}} {\|g(s, z)\|( {\phi^{\varepsilon ,\delta}  - 1} )\nu(dz)ds} }dt.
		\end{aligned}
		\end{eqnarray}
		Take similar manner in \eqref{lemma3.7}, by the definition of the stopping time, it implies that
		$$\mathbb{E}[{|{\mathcal{H}}_2|}] \le c_{18}\sqrt{\varepsilon h},\quad \mathbb{E}[{|{\mathcal{H}}_5|}] \le c_{19}\sqrt{\varepsilon}.$$
		Thus, by \cite[Lemma 3.4]{2013Large}, it can be concluded  from (\ref{lemma3.25}) that (\ref{lemma3.18}) holds.

Next, we define another stopping times
		{$\hat \tau_R^\varepsilon:=inf\{t\in[0,T]:| \hat{X}_{t}^{\varepsilon ,\delta} |+|  \tilde{X}_{t}^{\varepsilon ,\delta}|>R\}$}
for any $R, \varepsilon>0$. 
Then, we  are reducing to show that \begin{eqnarray}\label{lemma3.73}
\mathbb{E}[\mathop {\sup }\limits_{t \in [ {0, T\wedge\hat\tau_R^\varepsilon} ]} {| \hat{X}_{t}^{\varepsilon ,\delta} - \tilde{X}_{t}^{\varepsilon ,\delta}|^2}] \le c\Delta,
\end{eqnarray}
where $c>0$ is a constant independent of $\varepsilon,\delta,\Delta$. 

Define that ${| \hat{X}_{t}^{\varepsilon ,\delta} - \tilde{X}_{t}^{\varepsilon ,\delta} |^2}=\mathcal{J}_t^{\varepsilon,\delta}$, according  to the  It\^o's formula,  it leads to that 
		\begin{eqnarray*}\label{lemma3.41}
			\mathcal{J}_t^{\varepsilon,\delta}:= {{\rm \mathcal{J}}_1} + {{\rm\mathcal{J}}_2} + {{\rm \mathcal{J}}_3} + {{\rm \mathcal{J}}_4} + {{\rm \mathcal{J}}_5} + {{\rm \mathcal{J}}_6} + {{\rm \mathcal{J}}_7},
		\end{eqnarray*}
	   and
	    \begin{eqnarray*}\label{lemma3.80}
	    \begin{aligned}
		 {{\rm \mathcal{J}}_1} &= 2\int_0^t {\langle {{\rm \mathcal{J}}_s^{\varepsilon ,\delta }} ,[ {b_1( {\hat{X}_{s}^{\varepsilon ,\delta},\mathcal{L}_{X^{\varepsilon ,\delta}_s},\hat{Y}_{s}^{\varepsilon ,\delta}} ) - b_1( {\hat{X}_{s( \Delta  )}^{\varepsilon ,\delta},\mathcal{L}_{X^{\varepsilon ,\delta}_{s(\Delta)}},\tilde y_s^{\varepsilon ,\delta }} )} ]\rangle ds}\\
		{{\rm \mathcal{J}}_2} &= 2\int_0^t {\langle  {{\rm \mathcal{J}}_s^{\varepsilon ,\delta }} ,[ {\sigma_{1}( \hat{X}_{s}^{\varepsilon ,\delta},\mathcal{L}_{X^{\varepsilon ,\delta}_s} ) - \sigma_{1}( { \tilde{X}_{s}^{\varepsilon ,\delta}},\mathcal{L}_{X^{\varepsilon ,\delta}_{s}} )} ]\psi_s^{\varepsilon ,\delta} \rangle ds} \\
    	{{\rm \mathcal{J}}_3}
    	&= 2\sqrt \varepsilon  \int_0^t {\langle  {{\rm \mathcal{J}}_s^{\varepsilon ,\delta }}, \sigma_{1}( \hat{X}_{s}^{\varepsilon ,\delta},\mathcal{L}_{X^{\varepsilon ,\delta}_s} ) dW_s\rangle}, \\
    	{{\rm \mathcal{J}}_4}&=  \varepsilon \int_0^t {|{\sigma_{1}}( \hat{X}_{s}^{\varepsilon ,\delta},\mathcal{L}_{X^{\varepsilon ,\delta}_s} )|^2ds},\\
        {{\rm \mathcal{J}}_5}
        &= \int_0^t {\int_{\mathbf{X}} {[ {( {{{\rm \mathcal{J}}_s^{\varepsilon ,\delta }}{\rm{ + }}\varepsilon g(s, \hat{X}_{s}^{\varepsilon ,\delta} ,\mathcal{L}_{X^{\varepsilon ,\delta}_s},z)} )^2} - ( {{\rm \mathcal{J}}_s^{\varepsilon ,\delta }} )^2 ]{N^{{{\phi_s^{\varepsilon ,\delta}} \mathord{\left/{\vphantom {{\phi^{\varepsilon ,\delta} \left( s \right)} \varepsilon }} \right.\kern-\nulldelimiterspace} \varepsilon }}}( {dzds} )} },\\
        {{\rm \mathcal{J}}_6}
        &= - \int_0^t {\int_{\mathbf{X}} {2\langle {{\rm \mathcal{J}}_s^{\varepsilon ,\delta }},g(s, { \tilde{X}_{s}^{\varepsilon ,\delta}} ,\mathcal{L}_{X^{\varepsilon ,\delta}_{s}},z)\rangle \phi_s^{\varepsilon ,\delta}\nu(dz)ds} }, \\
        {{\rm \mathcal{J}}_7}&= -\int_0^t {\int_{\mathbf{X}} {2\langle {{\rm \mathcal{J}}_s^{\varepsilon ,\delta }},[ {g(s, \hat{X}_{s}^{\varepsilon ,\delta},\mathcal{L}_{X^{\varepsilon ,\delta}_s},z ) - g(s,{ \tilde{X}_{s}^{\varepsilon ,\delta}},\mathcal{L}_{x^{\varepsilon ,\delta}_{s}},z )} ]\rangle \nu(dz)ds} } \\
        	\end{aligned}
        \end{eqnarray*}
With Assumption {(\textbf{A1})}, we can see that
		\begin{eqnarray}\label{lemma3.42}
		\begin{aligned}
		{{\rm \mathcal{J}}_1} 
		&\le \int_{0}^t ( {{\rm \mathcal{J}}_s^{\varepsilon ,\delta }} )^2ds  + {{\rm \mathcal{J}}_{11}},\\
		{{\rm \mathcal{J}}_2} 
		&\le \int_0^t ( {{\rm \mathcal{J}}_s^{\varepsilon ,\delta }} )^2( {\psi_s^{\varepsilon ,\delta}} )^2ds  + {C_1}\int_0^t ( {{\rm \mathcal{J}}_s^{\varepsilon ,\delta }} )^2ds,
		\end{aligned}
		\end{eqnarray}
where
		\[
			{{\rm \mathcal{J}}_{11}}
:=\int_0^t  | {b_1( {\hat{X}_{s}^{\varepsilon ,\delta},\mathcal{L}_{X^{\varepsilon ,\delta}_s},\hat{Y}_{s}^{\varepsilon ,\delta}} ) - b_1( {\hat{X}_{s( \Delta  )}^{\varepsilon ,\delta},\mathcal{L}_{X^{\varepsilon ,\delta}_{s(\Delta)}},\tilde Y_s^{\varepsilon ,\delta }} )}  |^2ds.
		\]
Rearrange the sum of ${{\rm \mathcal{J}}_5}$, ${{\rm \mathcal{J}}_6}$, and ${{\rm \mathcal{J}}_7}$  as follows
		\begin{eqnarray*}
		{{\rm \mathcal{J}}_5} + {{\rm \mathcal{J}}_6}+{{\rm \mathcal{J}}_7}:= {{\rm \mathcal{J}}_{51}} + {{\rm \mathcal{J}}_{52}} + {{\rm \mathcal{J}}_{53}},
        \end{eqnarray*}
    where
		\begin{eqnarray*}\label{lemma3.44}
		\begin{aligned}
		{{\rm \mathcal{J}}_{51}} 
		&= {\varepsilon ^2}\int_0^t {\int_{\mathbf{X}} {  {g^2( s,\hat{X}_{s}^{\varepsilon ,\delta},\mathcal{L}_{X^{\varepsilon ,\delta}_s} z )}  {N^{{{\phi_s^{\varepsilon ,\delta} } \mathord{\left/{\vphantom {{\phi^{\varepsilon ,\delta} } \varepsilon }} \right.\kern-\nulldelimiterspace} \varepsilon }}}\left( {dzds} \right)} },\\
		{{\rm \mathcal{J}}_{52}}
		&=\int_0^t {\int_{\mathbf{X}} [ {2\varepsilon {{\rm \mathcal{J}}_s^{\varepsilon ,\delta }}g(s, \hat{X}_{s}^{\varepsilon ,\delta} ,\mathcal{L}_{X^{\varepsilon ,\delta}_s}, z)} ]{\tilde N}^{{\phi_s^{\varepsilon ,\delta} } \mathord{\left/{\vphantom {{\phi^{\varepsilon ,\delta} } \varepsilon }} \right.\kern-\nulldelimiterspace} \varepsilon }( {dzds} ) },\\
		{{\rm \mathcal{J}}_{53}}&= \int_0^t {\int_{\mathbf{X}} {2\langle {{\rm \mathcal{J}}_s^{\varepsilon ,\delta }},[ {g(s, \hat{X}_{s}^{\varepsilon ,\delta},\mathcal{L}_{X^{\varepsilon ,\delta}_s} , z) - g(s,{ \tilde{X}_{s}^{\varepsilon ,\delta}},\mathcal{L}_{x^{\varepsilon ,\delta}_{{s}}} , z)} ]\rangle ( {\phi_s^{\varepsilon ,\delta}  - 1} )\nu(dz)ds} }.
		\end{aligned}
		\end{eqnarray*}
It implies form Assumption {(\textbf{A1})} that
		\begin{eqnarray}\label{lemma3.45}
		{{\rm \mathcal{J}}_{53}} &\le& \int_0^t {\int_{\mathbf{X}} ( {{\rm \mathcal{J}}_s^{\varepsilon ,\delta }} )^2| \phi_s^{\varepsilon ,\delta} - 1|\nu(dz)ds } +\int_0^t {\int_{\mathbf{X}} ( {{\rm \mathcal{J}}_s^{\varepsilon ,\delta }} )^2{\|g(s, z)\|}| \phi_s^{\varepsilon ,\delta} - 1|\nu(dz)ds }.
		\end{eqnarray}
		With aid of  the  Gronwall's inequality, it deduce from (\ref{lemma3.42}) to (\ref{lemma3.45}), 
		\begin{eqnarray*}\label{lemma3.46}
			{| \hat{X}_{t}^{\varepsilon ,\delta} - \tilde{X}_{t}^{\varepsilon ,\delta} |^2} \le {e^{\left(   2t + {C_1}{M_\psi } + M_\phi \right)}}\left\{ {{{\rm \mathcal{J}}_{11}} + {{\rm \mathcal{J}}_3} + {{\rm \mathcal{J}}_4} + {{\rm \mathcal{J}}_{51}} + {{\rm \mathcal{J}}_{52}} } \right\},
		\end{eqnarray*}
which leads to that		
		\begin{eqnarray}\label{lemma3.47}
		\mathbb{E}[\mathop {\sup }\limits_{t \in [ {0, T\wedge\hat\tau_R^\varepsilon} ]} {|  \hat{X}_{t}^{\varepsilon ,\delta} - \tilde{X}_{t}^{\varepsilon ,\delta}  |^2}] \le c_{20}\mathbb{E}[\mathop {\sup }\limits_{t \in [ {0, T\wedge\hat\tau_R^\varepsilon} ]} \left( {{{\rm \mathcal{J}}_{11}} + {{\rm \mathcal{J}}_3} + {{\rm \mathcal{J}}_4} + {{\rm \mathcal{J}}_{51}} + {{\rm \mathcal{J}}_{52}}} \right)],
		\end{eqnarray}
where we  choose the constant $c_{20}\ge{e^{\left(   2t + {C_1}{M_\psi } +M_\phi \right)}}$.
According to the  \cite[Lemma 5.8]{Hong2021Central} that
		for $\varepsilon, \delta> 0$,   $t \in [0, T\wedge\tau_R^\varepsilon]$, and $(\psi^{\varepsilon ,\delta}, \phi^{\varepsilon ,\delta})\in \mathcal{U}^{m}$,
		\begin{eqnarray}\label{lemma3.86}
			\mathbb{E}[\int_0^{ T\wedge\tau_R^\varepsilon}{| {{ \hat{Y}_{t}^{\varepsilon ,\delta}} - \tilde Y_t^{\varepsilon ,\delta }} |^2}dt ]\le c_{21}\Delta.
		\end{eqnarray} 		
Then on account of Assumption {(\textbf{A1})}, \lemref{lem2}, estimates (\ref{lemma3.18}) and (\ref{lemma3.86}), it follows that
		\begin{eqnarray}\label{lemma3.48}
		\begin{aligned}
		\mathbb{E}[\mathop {\sup }\limits_{t \in [ {0, T\wedge\hat\tau_R^\varepsilon} ]}  {{{\rm \mathcal{J}}_{11}}}]  &\le\mathop {\sup }\limits_{t \in [ {0, T\wedge\hat\tau_R^\varepsilon} ]} \int_0^t [ {b_1( {\hat{X}_{s}^{\varepsilon ,\delta},\mathcal{L}_{X^{\varepsilon ,\delta}_s},\hat{Y}_{s}^{\varepsilon ,\delta}} ) - b_1( {\hat{X}_{s\left( \Delta  \right)}^{\varepsilon ,\delta},\mathcal{L}_{X^{\varepsilon ,\delta}_{s(\Delta)}},\tilde Y_s^{\varepsilon ,\delta }} )} ]^2ds \\
		&\le {C_1}\mathbb{E} \int_0^{\hat\tau_R^\varepsilon} \big[ ( \hat{X}_{s}^{\varepsilon ,\delta} - \hat{X}_{s\left( \Delta  \right)}^{\varepsilon ,\delta} )^2 + (\hat{Y}_{s}^{\varepsilon ,\delta} - \tilde Y_s^{\varepsilon ,\delta } ) ^2 +\mathbb{W}^2_2(\mathcal{L}_{X^{\varepsilon ,\delta}_s},\mathcal{L}_{X^{\varepsilon ,\delta}_{s(\Delta)}})\big]ds \\
		&\le c_{22}\Delta.
		\end{aligned}
		\end{eqnarray}
		By Assumptions {(\textbf{A1})} and (\textbf{A4}), we can see that
		\begin{eqnarray}\label{lemma3.50}
		\begin{aligned}
		\mathbb{E}[\mathop {\sup }\limits_{t \in [ {0, T\wedge\hat\tau_R^\varepsilon} ]}  {{{\rm \mathcal{J}}_4}}]  & \le \varepsilon {C_2 T},\\
		\mathbb{E}[\mathop {\sup }\limits_{t \in [ {0, T\wedge\hat\tau_R^\varepsilon} ]}  {{{\rm \mathcal{J}}_{51}}} ] &\le \mathbb{E}\Big[\mathop {\sup }\limits_{t \in [ {0, T\wedge\hat\tau_R^\varepsilon} ]}  \varepsilon^2 \int_0^t \int_{\mathbf{X}}  g^2(s, \hat{X}_{s}^{\varepsilon ,\delta},\mathcal{L}_{X^{\varepsilon ,\delta}_s}, z )  {N^{{{\phi_s^{\varepsilon ,\delta} } \mathord{\left/{\vphantom {{\phi_s^{\varepsilon ,\delta} } \varepsilon }} \right.\kern-\nulldelimiterspace} \varepsilon }}}\left( {dzds} \right)  \Big]\\
		&\le \varepsilon \mathbb{E}\Big[ \int_0^{ T\wedge\hat\tau_R^\varepsilon} \int_{\mathbf{X}} ( 1 + | \hat{X}_{s}^{\varepsilon ,\delta} |^2 +\mathcal{L}_{X^{\varepsilon ,\delta}_s}(|\cdot|^2))\|g(s, z)\|\phi_s^{\varepsilon ,\delta} \nu(dz)ds   \Big] \\
		&\le \varepsilon c {M_\phi} + \varepsilon R^2{M_\phi}.
		\end{aligned}
		\end{eqnarray}		
With aid of  the Burkholder-Davis-Gundy inequality, it leads to 
		\begin{eqnarray}\label{lemma3.49}
		\begin{aligned}
		\mathbb{E}[\mathop {\sup }\limits_{t \in [ {0, T\wedge\hat\tau_R^\varepsilon} ]} {{{\rm \mathcal{J}}_3}}  ]&\le 8\mathbb{E}{\Big\{ {\varepsilon \int_0^{ T\wedge\hat\tau_R^\varepsilon} {{( {{\rm \mathcal{J}}_s^{\varepsilon ,\delta }} )^2}|{\sigma_{1}}( \hat{X}_{s}^{\varepsilon ,\delta} ,\mathcal{L}_{X^{\varepsilon ,\delta}_s})|^2ds} } \Big\}^{{1 \mathord{\left/{\vphantom {1 2}} \right.\kern-\nulldelimiterspace} 2}}}\\
		&\le \frac{1}{{8c_{20}}}\mathbb{E}[\mathop {\sup }\limits_{t \in [ {0, T\wedge\hat\tau_R^\varepsilon} ]} ( {{\rm \mathcal{J}}_t^{\varepsilon ,\delta }} )^2] + 128\varepsilon c_{23}(1+R^2+\mathbb{E}[\sup_{t\in[0,T]}|x_t^{\varepsilon ,\delta}|^2]),\\
			\mathbb{E}[\mathop {\sup }\limits_{t \in[ {0, T\wedge\hat\tau_R^\varepsilon} ]}  {{{\rm \mathcal{J}}_{52}}} ] &\le 4\mathbb{E}\left[ {{\rm \mathcal{J}}_{52}^2} \right]_{T\wedge\hat\tau_R^\varepsilon}^{{1 \mathord{\left/{\vphantom {1 2}} \right.\kern-\nulldelimiterspace} 2}}\\
			&\le 8\mathbb{E}{\Big[ {\varepsilon \mathop {\sup }\limits_{t \in [ {0, T\wedge\hat\tau_R^\varepsilon} ]} | {{\rm \mathcal{J}}_t^{\varepsilon ,\delta }} |^2\int_0^{ T\wedge\hat\tau_R^\varepsilon} {\int_{\mathbf{X}} {{g^2}(s, \hat{X}_{s}^{\varepsilon ,\delta},\mathcal{L}_{X^{\varepsilon ,\delta}_s}, z )\phi_s^{\varepsilon ,\delta} \nu( dz )ds} } } \Big]^{{1 \mathord{\left/{\vphantom {1 2}} \right.\kern-\nulldelimiterspace} 2}}}\\
			&\le \frac{1}{{8c_{20}}}\mathbb{E}[\mathop {\sup }\limits_{t \in [ {0, T\wedge\hat\tau_R^\varepsilon} ]} {| {{\rm \mathcal{J}}_t^{\varepsilon ,\delta }}|^2}] +128 c_{23}\varepsilon M_\phi\mathbb{E}[\sup_{t\in[0,T]}|X_t^{\varepsilon ,\delta}|^2]\\
			&\quad +128 c_{23}\varepsilon M_\phi R^2+ 128c_{23}\varepsilon  M_\phi T.
			\end{aligned}
		\end{eqnarray}				
Then by estimates (\ref{lemma3.47}),  (\ref{lemma3.48})--(\ref{lemma3.49}), we can conclude that
		\begin{eqnarray*}
		\mathbb{E}[\mathop {\sup }\limits_{t \in [ {0, T\wedge\hat\tau_R^\varepsilon} ]} {| \hat{X}_{t}^{\varepsilon ,\delta} - \tilde{X}_{t}^{\varepsilon ,\delta}|^2}] \le c_{24}\Delta  + \frac{1}{4}\mathbb{E}[\mathop {\sup }\limits_{t \in [ {0, T\wedge\hat\tau_R^\varepsilon} ]} {| \hat{X}_{t}^{\varepsilon ,\delta} - \tilde{X}_{t}^{\varepsilon ,\delta} |^2}],
		\end{eqnarray*}
		which leads to that
				\begin{eqnarray}\label{lemma3.53}
		\mathbb{E}[\mathop {\sup }\limits_{t \in [ {0, T\wedge\hat\tau_R^\varepsilon} ]} {| \hat{X}_{t}^{\varepsilon ,\delta} - \tilde{X}_{t}^{\varepsilon ,\delta}|^2}] \le \frac{4}{3}c_{24}\Delta,
		\end{eqnarray}
		therefore, the estimate \eqref{lemma3.73} is obtained.
		
	Then we only need to show that for $t \in [0, T]$ and $\varepsilon, \delta>0$ small enough,  $(\psi^{\varepsilon ,\delta}, \phi^{\varepsilon ,\delta})\in \mathcal{U}^{m}$,   
				\begin{eqnarray}\label{lemma3.55}
							\mathbb{E}[\mathop {\sup }\limits_{t \in [ {0, T\wedge\hat\tau_R^\varepsilon} ]} {|\tilde{X}_{t}^{\varepsilon ,\delta} - \hat X_t^{ \varepsilon}|^2}] \le \frac{4}{3}c_{24}\Delta.
				\end{eqnarray}
				Define that $\tilde{X}_{t}^{\varepsilon ,\delta}-\hat X_t^{ \varepsilon} =  {\mathcal{K}}^{\varepsilon ,\delta }_t$, 
			{\begin{eqnarray*}\label{lemma3.56}
					{{\rm \mathcal{K}}_t^{\varepsilon ,\delta }}  &=:& {{\rm \mathcal{K}}_1} +{{\rm \mathcal{K}}_2}+ {{\rm \mathcal{K}}_3} + {{\rm \mathcal{K}}_4},
			\end{eqnarray*}}	
				where
				\begin{eqnarray*}\label{lemma3.57}
				\begin{aligned}
					{{\rm \mathcal{K}}_1}&=\int_0^t [{b_1( { \hat{X}_{s(\Delta)}^{\varepsilon ,\delta}},\mathcal{L}_{x^{\varepsilon ,\delta}_{s(\Delta)}},\tilde Y_s^{\varepsilon ,\delta } ) - \bar b_1( { \hat{X}_{s(\Delta)}^{\varepsilon ,\delta}},\mathcal{L}_{X^{\varepsilon ,\delta}_{s(\Delta)}})}]ds,\cr
					{{\rm \mathcal{K}}_2}&=\int_0^t [\bar b_1( { \hat{X}_{s(\Delta)}^{\varepsilon ,\delta}} ,\mathcal{L}_{X^{\varepsilon ,\delta}_{s(\Delta)}}) - \bar b_1( {\hat X_s^{ \varepsilon}} ,\mathcal{L}_{\bar X_s}) ]ds,\cr
					{{\rm \mathcal{K}}_3}&=\int_0^t [\sigma_{1}( \tilde{X}_{s}^{\varepsilon ,\delta},\mathcal{L}_{X^{\varepsilon ,\delta}_{{s}}} ) - \sigma_{1}( \hat X_s^{ \varepsilon},\mathcal{L}_{\bar X_s} )]\psi_s^{\varepsilon ,\delta}ds,\cr
					{{\rm \mathcal{K}}_4}&= \int_{[0, T]\times \mathbf{X}} {\big[ {g(s, {\tilde{X}_{s}^{\varepsilon ,\delta}} ,\mathcal{L}_{X^{\varepsilon ,\delta}_{s}}, z) - g(s,{\hat X_s^{ \varepsilon}},\mathcal{L}_{\bar X_s}, z)} \big]( {\phi_s^{\varepsilon ,\delta}  - 1} )\nu(dz)ds},
			     \end{aligned}
				\end{eqnarray*}
				
				According to Assumption {(\textbf{A1})}, \lemref{lem3} and  estimate (\ref{lemma3.18}), we have
				\begin{eqnarray*}\label{lemma3.58}
					{\left| {{{\rm \mathcal{K}}_2}} \right|^2} &\le& {\Big\{ {\int_0^t {\big(\bar b_1( { \hat{X}_{s(\Delta)}^{\varepsilon ,\delta}},\mathcal{L}_{X^{\varepsilon ,\delta}_{s(\Delta)}} ) - \bar b_1( {\hat X_s^{ \varepsilon}},\mathcal{L}_{\bar X_s} ) \big)ds} } \Big\}^2}\cr
					&\le& {\Big\{ {\int_0^t {\big(\bar b_1( \tilde{X}_{s}^{\varepsilon ,\delta},\mathcal{L}_{X^{\varepsilon ,\delta}_{s(\Delta)}} ) - \bar b_1( {\hat X_s^{ \varepsilon}},\mathcal{L}_{X^{\varepsilon ,\delta}_{s}} ) \big)ds} } \Big\}^2}+O(\Delta)\cr
					&\le& {C_1} {\int_0^t {{( {{\rm \mathcal{K}}_s^{\varepsilon ,\delta }} )^2}ds} }+O(\Delta) ,
				\end{eqnarray*}
				likewise, it deduces that 
				\begin{eqnarray*}
					{\left| {{{\rm \mathcal{K}}_3}} \right|^2} &\le& {\Big\{ \int_0^t \big[ {\sigma_{1}( {\tilde{X}_{s}^{\varepsilon ,\delta}} ,\mathcal{L}_{X^{\varepsilon ,\delta}_{{s}}}) - \sigma_{1}( {\hat X_s^{ \varepsilon}} ,\mathcal{L}_{\bar X_s})}\big] \psi_s^{\varepsilon ,\delta}ds  \Big\}^2}\cr
					&\le& {C_1} {\int_0^t {{( {{\rm \mathcal{K}}_s^{\varepsilon ,\delta }} )^2|\psi_s^{\varepsilon ,\delta}|^2}ds} }+\Delta.
				\end{eqnarray*}
				Furthermore, by the Assumption (\textbf{A1}), it has
				\begin{eqnarray*}\label{lemma3.59}
					{\left| {{{\rm \mathcal{K}}_4}} \right|^2} &\le& 2\int_0^t {\int_{\mathbf{X}} { {{\rm \mathcal{K}}_4}\big[ {g(s,  \tilde{X}_{s}^{\varepsilon ,\delta} ,\mathcal{L}_{X^{\varepsilon ,\delta}_{{s}}}, z) - g(s, \hat X_s^{ \varepsilon},\mathcal{L}_{\bar X_s}, z )} \big]( {\phi_s^{\varepsilon ,\delta}  - 1} )\nu(dz)ds} } \cr
					&\le& \int_0^t {\int_{\mathbf{X}} { {| {{{\rm \mathcal{K}}_4}} |^2}| {\phi_s^{\varepsilon ,\delta}  - 1} |\nu(dz)ds} } \cr
					&& + \int_0^t {\int_{\mathbf{X}} { {| {{\rm \mathcal{K}}_s^{\varepsilon ,\delta }} |^2}\|g(s, z)\|| {\phi_s^{\varepsilon ,\delta}  - 1} |\nu(dz)ds} } \cr
					&\le& {e^{M_\phi}}\int_0^t {\int_{\mathbf{X}} { {| {{\rm \mathcal{K}}_s^{\varepsilon ,\delta }} |^2}\|g(s, z)\|| {\phi_s^{\varepsilon ,\delta}  - 1} |\nu(dz)ds} } .
				\end{eqnarray*}
				Then by the Grownwall lemma,
				\begin{eqnarray*}\label{lemma3.60}
					{| {\mathcal{K}}^{\varepsilon ,\delta }_t |^2} \le {e^{( {{C_1}{M_\psi }T +C_1 T+ {e^{M_\phi}}M_\phi} )}} ({\rm \mathcal{K}}_1)^2+{e^{( {{C_1}{M_\psi }T +C_1 T+ {e^{M_\phi}}M_\phi} )}}O(\Delta) .
				\end{eqnarray*}
				Set $c_{25}\ge{e^{( {{C_1}{M_\psi }T +C_1 T+ {e^{M_\phi}}M_\phi} )}}$, then,
				\begin{eqnarray*}\label{lemma3.61}
					\mathbb{E}[\mathop {\sup }\limits_{t \in \left[ {0, T\wedge\hat\tau_R^\varepsilon} \right]} {| {\mathcal{K}}^{\varepsilon ,\delta }_t |^2} ]\le c_{25}\mathbb{E}\mathop {\sup }\limits_{t \in \left[ {0,T} \right]}  ({\rm \mathcal{K}}_1)^2 +c_{25}O(\Delta).
				\end{eqnarray*}
				Construct ${{\tilde y}^{{ \hat{X}_{k\Delta}^{\varepsilon ,\delta}},\tilde Y_{k\Delta}^{\varepsilon ,\delta }}}\left( {\frac{s}{\delta }} \right)$ as follows,
				\begin{eqnarray*}\label{lemma3.63}
					{{\tilde Y}^{{ \hat{X}_{k\Delta}^{\varepsilon ,\delta}},\tilde Y_{k\Delta}^{\varepsilon ,\delta }}}\left( {\frac{s}{\delta }} \right) &=& {{\tilde Y}^{\varepsilon ,\delta } }_{k\Delta }  + \int_0^{\frac{s}{\delta}} b_2(  \hat{X}_{k\Delta}^{\varepsilon ,\delta},\mathcal{L}_{X^{\varepsilon ,\delta}_{{k\Delta}}},\tilde Y_u^{ \hat{X}_{k\Delta}^{\varepsilon ,\delta},\tilde Y_{k\Delta}^{\varepsilon ,\delta }} )du+ \int_0^{\frac{s}{\delta}} {\sigma_{2}( { \hat{X}_{k\Delta}^{\varepsilon ,\delta}},\mathcal{L}_{X^{\varepsilon ,\delta}_{k\Delta}},\tilde Y_u^{ \hat{X}_{k\Delta}^{\varepsilon ,\delta},\tilde Y_{k\Delta}^{\varepsilon ,\delta }})dW_u},
				\end{eqnarray*}								for $0 \leqslant k\leqslant[\frac{t}{\Delta}]-1$.
				Furthermore, we have that $\mathop {\sup }\limits_{0 \le t \le T} {({\rm \mathcal{K}}_1)^2} \le 	{\mathcal I_{11}}+{\mathcal I_{12}},$
				where
				\begin{eqnarray}\label{lemma3.64}
				{\mathcal I_{11}}=8\mathbb{E}\mathop {\sup }\limits_{0 \le t \le T} {\Big[ {\sum\limits_{k = 0}^{\left[ {\frac{t}{\Delta}} \right]-1} {\int_{k\Delta }^{\left( {k + 1} \right)\Delta } { ( {b_1( { \hat{X}_{k\Delta}^{\varepsilon ,\delta}},\mathcal{L}_{X^{\varepsilon ,\delta}_{{k\Delta}}},\tilde Y_s^{\varepsilon ,\delta }) - \bar b_1( { \hat{X}_{k\Delta}^{\varepsilon ,\delta}},\mathcal{L}_{X^{\varepsilon ,\delta}_{k\Delta}} )} )ds} } } \Big]^2},
				\end{eqnarray}		
				and
				\begin{eqnarray*}
					{\mathcal I_{12}}=8\mathbb{E}\mathop {\sup }\limits_{0 \le t \le T} {\Big[  {\int_{\left[ {\frac{t}{\Delta}} \right]\Delta }^{{t} } { ( {b_1( { \hat{X}_{k\Delta}^{\varepsilon ,\delta}},\mathcal{L}_{x^{\varepsilon ,\delta}_{k\Delta}},\tilde Y_s^{\varepsilon ,\delta } ) - \bar b_1({ \hat{X}_{k\Delta}^{\varepsilon ,\delta}} ,\mathcal{L}_{X^{\varepsilon ,\delta}_{k\Delta}})} )ds} }  \Big]^2}.
				\end{eqnarray*}	
				Then change the time scale, it deduces that
				\begin{eqnarray}\label{lemma3.65}
				{{ \mathcal{K}}_{11}}
				&\le& 8{\delta ^2}{[ {\frac{T}{\Delta}} ]^2}\mathop {\sup }\limits_{0 \le k \le [ {\frac{T}{\Delta}} ]-1} \mathbb{E}{\Big[\Big| {\int_0^{\frac{\Delta}{\delta}} {\big[ {b_1(  \hat{X}_{k\Delta}^{\varepsilon ,\delta},\mathcal{L}_{X^{\varepsilon ,\delta}_{k\Delta}},{\tilde Y}_s^{{ \hat{X}_{k\Delta}^{\varepsilon ,\delta}},\tilde Y_{k\Delta}^{\varepsilon ,\delta }} ) - \bar b_1( { \hat{X}_{k\Delta}^{\varepsilon ,\delta}} ,\mathcal{L}_{X^{\varepsilon ,\delta}_{k\Delta}})} \big]ds} } \Big|^2\Big]}\cr
				&\le& 8{\delta ^2}{\left[ {\frac{T}{\Delta}} \right]^2}\mathop {\max }\limits_{0 \le k \le \left[ {\frac{T}{\Delta}} \right]-1} {\rm \mathcal{K}}_k^\delta.
				\end{eqnarray}									
			By	the	Cauchy-Schwarz's inequality and  exponential ergodicity of the fast component  \cite[Proposition 3.7]{M2021Strong}, it  leads to that
				\begin{eqnarray}\label{lemma3.66}
				{\rm \mathcal{K}}_k^\delta &=& \int_0^{\frac{\Delta}{\delta}} \int_\tau ^{\frac{\Delta}{\delta}} \mathbb{E}\Big\{ \big[ {b_1(  \hat{X}_{k\Delta}^{\varepsilon ,\delta},\mathcal{L}_{X^{\varepsilon ,\delta}_{k\Delta}},{\tilde Y}_s^{{ \hat{X}_{k\Delta}^{\varepsilon ,\delta}},\tilde Y_{k\Delta}^{\varepsilon ,\delta }} ) - \bar b_1( { \hat{X}_{k\Delta}^{\varepsilon ,\delta}} ,\mathcal{L}_{X^{\varepsilon ,\delta}_{k\Delta}})} \big]\Big.\cr
				&&\quad\Big.\big[ b_1(  \hat{X}_{k\Delta}^{\varepsilon ,\delta},\mathcal{L}_{X^{\varepsilon ,\delta}_{k\Delta}},{\tilde Y}_\tau^{{ \hat{X}_{k\Delta}^{\varepsilon ,\delta}},\tilde Y_{k\Delta}^{\varepsilon ,\delta }} ) - \bar b_1( { \hat{X}_{k\Delta}^{\varepsilon ,\delta}} ,\mathcal{L}_{X^{\varepsilon ,\delta}_{k\Delta}}) \big] \Big\}  dsd\tau \cr
				&=&\int_0^{\frac{\Delta}{\delta}} \int_\tau ^{\frac{\Delta}{\delta}} {\mathbb{E}^y}\Big\{ \big[ {b_1\big( {{ \hat{X}_{k\Delta}^{\varepsilon ,\delta}},\mathcal{L}_{X^{\varepsilon ,\delta}_{k\Delta}},{\tilde y}_\tau^{{ \hat{X}_{k\Delta}^{\varepsilon ,\delta}},\tilde Y_{k\Delta}^{\varepsilon ,\delta }}} \big) - \bar b_1\big( { \hat{X}_{k\Delta}^{\varepsilon ,\delta}} ,\mathcal{L}_{X^{\varepsilon ,\delta}_{k\Delta}}\big)} \big]\Big.\cr
				&&\quad\Big.{\mathbb{E}^{{Y^{X,Y}}\left( \tau  \right)}} \big[ b_1\big( { \hat{X}_{k\Delta}^{\varepsilon ,\delta}},\mathcal{L}_{X^{\varepsilon ,\delta}_{k\Delta}},{\tilde Y}_{s-\tau}^{{ \hat{X}_{k\Delta}^{\varepsilon ,\delta}},\tilde Y_{k\Delta}^{\varepsilon ,\delta }} \big) - \bar b_1\big( { \hat{X}_{k\Delta}^{\varepsilon ,\delta}} ,\mathcal{L}_{X^{\varepsilon ,\delta}_{k\Delta}}\big) \big] \Big\}dsd\tau   \cr
				&\le& \int_0^{\frac{\Delta}{\delta}} \int_\tau ^{\frac{\Delta}{\delta}} {\Big\{ {{\mathbb{E}^Y}{ \big[ {b_1( { \hat{X}_{k\Delta}^{\varepsilon ,\delta}},\mathcal{L}_{X^{\varepsilon ,\delta}_{k\Delta}},{\tilde Y}_\tau^{{ \hat{X}_{k\Delta}^{\varepsilon ,\delta}},\tilde Y_{k\Delta}^{\varepsilon ,\delta }} ) - \bar b_1( { \hat{X}_{k\Delta}^{\varepsilon ,\delta}} ,\mathcal{L}_{X^{\varepsilon ,\delta}_{k\Delta}})} \big]^2}} \Big\}^{{1 \mathord{\left/{\vphantom {1 2}} \right.\kern-\nulldelimiterspace} 2}}}\cr
				&&\quad{\Big\{ {{\mathbb{E}^Y}{\big\{ {{\mathbb{E}^{{Y^{X,Y}}\left( \tau  \right)}} {\big[ {b_1\big( { \hat{X}_{k\Delta}^{\varepsilon ,\delta}},\mathcal{L}_{X^{\varepsilon ,\delta}_{k\Delta}},{\tilde Y}_{s-\tau}^{{ \hat{X}_{k\Delta}^{\varepsilon ,\delta}},\tilde Y_{k\Delta}^{\varepsilon ,\delta }} \big) - \bar b_1\big( { \hat{X}_{k\Delta}^{\varepsilon ,\delta}} ,\mathcal{L}_{X^{\varepsilon ,\delta}_{k\Delta}}\big)} \big]} } \big\}^2}} \Big\}^{{1 \mathord{\left/{\vphantom {1 2}} \right.\kern-\nulldelimiterspace} 2}}}  dsd\tau \cr
				&\le& c_{26}\int_0^{\frac{\Delta}{\delta}} {\int_\tau ^{\frac{\Delta}{\delta}} {{e^{ - \frac{{c\left( {s - \tau } \right)}}{2}}}dsd\tau } } \cr
				&\le& c_{26}\big( {\frac{4}{{{\eta ^2}}}{e^{ - \frac{{c \Delta }}{{2\delta }}}}-\frac{4}{{{\eta ^2}}} + \frac{{2\Delta }}{{\eta \delta }}} \big).
				\end{eqnarray}	
				By Assumption (\textbf{A1}), the definition of stopping times, and  (\ref{lemma3.15})
				\begin{eqnarray}\label{lemma3.67}
				{{\rm \mathcal{K}}_{12}} &\le&  c_{27}{\Delta }.
				\end{eqnarray}	 							
				Hence, from (\ref{lemma3.64}) to (\ref{lemma3.67}),
				\begin{eqnarray*}\label{lemma3.68}
		({\rm \mathcal{K}}_1)^2  \le c_{26}\frac{\delta }{\Delta }+c_{27}{\Delta }.
				\end{eqnarray*}
			Hence, we have 
			\begin{eqnarray}\label{lemma3.69}
			\mathbb{E}[\mathop {\sup }\limits_{t \in \left[ {0, T\wedge\hat\tau_R^\varepsilon} \right]} {| {\mathcal{K}}^{\varepsilon ,\delta }_t |^2} ]\le c_{25}\mathbb{E}\mathop {\sup }\limits_{t \in \left[ {0,T} \right]} ({\rm \mathcal{K}}_1)^2+c_{25}O(\Delta) \le c_{28}\frac{\delta }{\Delta }+c_{28}O(\Delta).
			\end{eqnarray}		
			{
				Combine \eqref{lemma3.53} and \eqref{lemma3.69}, we have 
				\begin{eqnarray}\label{lemma3.71}
				\mathbb{E}[\mathop {\sup }\limits_{t \in \left[ {0, T\wedge\hat\tau_R^\varepsilon} \right]} {|\hat{X}_{t}^{\varepsilon ,\delta}-\hat X_t^{\varepsilon}|^2} ]\le c(\frac{\delta }{\Delta }+\Delta),
				\end{eqnarray}
			where $c$ is independent of $\varepsilon,\delta,\Delta$.}	
					
		{According to  the definition of stopping times $\hat\tau_R^\varepsilon$, we can get for any $r>0$
			\begin{eqnarray}\label{lemma3.70}
			\mathbb{P}\big(\sup _{t \in[0, T]}|\hat{X}_{t}^{\varepsilon ,\delta}-\hat X_t^{\varepsilon}| \geqslant r\big) 
			&\le & \mathbb{P}\left(T>\hat{\tau}_R^{\varepsilon}\right)+\mathbb{P}\big(\sup _{t \in[0, T]}|\hat{X}_{t}^{\varepsilon ,\delta}-\hat X_t^{\varepsilon}| \geqslant r, T \le \hat{\tau}_R^{\varepsilon}\big) \cr
			&\le & \mathbb{P}\big(\mathop {\sup }\limits_{t \in \left[ {0,T} \right]}| \hat{X}_{t}^{\varepsilon ,\delta} |+\mathop {\sup }\limits_{t \in \left[ {0,T} \right]}|  \tilde{X}_{t}^{\varepsilon ,\delta}|>R\big) \cr
			&&+\mathbb{P}\big(\sup _{t \in\left[0, T \wedge \hat\tau_R^{\varepsilon}\right]}|\hat{X}_{t}^{\varepsilon ,\delta}-\hat X_t^{\varepsilon}| \le r\big).
		\end{eqnarray}
Firstly, for any fixed $R>0$, 	by estimates \eqref{lemma3.102} and \eqref{3-2}, the second part could  be small enough by choosing suitable $\Delta$ such as $\Delta=\sqrt{\delta}$, so that $\frac{\delta}{\Delta}$ small enough. Next, we will let  $R\to \infty$.}


{According to the Step 1, if  $(\psi^{\varepsilon ,\delta}, \phi^{\varepsilon ,\delta})\in \mathcal{U}^{m}$ such that $(\psi^{\varepsilon ,\delta}, \phi^{\varepsilon ,\delta})$ weakly converges  to $(\psi, \phi)$ as $\varepsilon \to 0$, then,  $\hat X^{\varepsilon}=\mathcal{G}^0 (\psi^{\varepsilon ,\delta}, \phi^{\varepsilon ,\delta})$ weakly converges to $\hat X=\mathcal{G}^0 (\psi, \phi)$ in $ \mathbf{D}$ as $\varepsilon \to 0$. 
	Then,  for any bounded continuous functions $h:\mathbf{D} \to\mathbb{R}$, we see that as $\varepsilon \to 0$
	\begin{eqnarray*}\label{3-82}
		|\mathbb{E}[h(\hat X^{\varepsilon ,\delta})]-\mathbb{E}[h(\hat X)]|&\le& |\mathbb{E}[h(\hat X^{\varepsilon ,\delta})]-\mathbb{E}[h(\hat X^{\varepsilon})]|+|\mathbb{E}[h(\hat X^{\varepsilon})]-\mathbb{E}[h(\hat X)]|\to 0,
	\end{eqnarray*} 
where we use the Portemanteau's theorem \cite[Theorem 13.16]{2020Klenke}.
	 Thus, we have obatined (\ref{step3}). }

	{\textbf{Step 3}.	With Step 1 and Step 2, it deduces from \cite[Theorem 4.4]{W2022Large} that $X^{\varepsilon ,\delta}$ satisfies  a large deviation principle on $\mathbf{D}$ with the good rate function $I: \mathbf{D}\rightarrow [0, \infty)$ defined in (\ref{rate}).}
	
	This proof is completed.\qed

	\section*{Statement of Contribution}
	Our work gives large deviations for a two-time scale McKean-Vlasov system with jumps. Different from previous general stochastic system, this McKean-Vlasov system does not only depends on the microcosmic location but also depends on the macrocosmic distribution. The novelty in this work is to treat this dependence. This large deviation result could provide theoretical framework for the long-time behavior for two-time scale McKean-Vlasov system in the real world.  
	
	\section*{Acknowledgments}
     This work was partly supported by  the NSF of China (Grant 12120101002), the NSF of China (Grant 12072264), the Fundamental Research Funds for the Central Universities, the Research Funds for Interdisciplinary Subject of Northwestern Polytechnical University,  the Shaanxi Provincial Key R\&D Program (Grants 2020KW-013, 2019TD-010).
	

\end{document}